# Integrating Wildfire Risk into Robust Forest Management Optimization

**Clemence Labarre (1,2), Denis Loustau (1), Jean-Christophe Domec (1), Mengistie Kindu (3), Thomas Knoke (3)**

- (1) Bordeaux Sciences-Agro, INRAE, UMR ISPA, Villenave d'Ornon, 33140, France.
- (2) Permanent address : IGN - Institut National de l'Information Géographique et Forestière, Saint-Médard en Jalles, 33160, France
- (3) Institute of Forest Management, TUM School of Life Sciences Weihenstephan, Life Science Systems, Technische Universität München, Hans-Carl-von-Carlowitz-Platz 2, 85354, Freising, Germany

  Clémence Labarre* (corresponding author) clemence.labarre@ign.fr

**Acknowledgements:**

French Ministry of Agriculture and INRAE
Nouvelle Aquitaine projects GRIFON and IMPACTS.

## Abstract

Increasingly severe and unpredictable extreme weather events—such as high-intensity forest fires—threaten long-term ecosystem service provision, suggesting a critical reassessment of current forest management planning. While planning frameworks exist, few incorporate optimization methods that account for the uncertainty and variability introduced by such extremes. Meanwhile, process-based models have advanced our understanding of fire impacts on forest dynamics and ecosystem services, yet their potential to inform strategic management decisions remains largely untapped. Focusing on intensive plantation forests, we propose a combined approach that uses fire simulation within a growth model to optimize forested landscape management under increasing fire risk. We demonstrate that incorporating fire behavior and its uncertainties significantly alters management outcomes, favoring diversified portfolios that balance high-yield but fire-prone regimes with more fire-resistant alternatives. While diversification improves robustness against fire for objectives like timber provision, it does not guarantee complementary benefits for carbon storage, even if carbon storage was considered as an objective. Our method exposes the limitations of current management regimes and identifies key areas for improved planning to better mitigate fire risk and strengthen long-term resilience.

## Highlights

- Coupled fire–growth modeling shows management can buffer timber losses.
- Robust, fire-aware optimization guides forest planning in complex systems.
- Diversified regimes secure timber, while carbon remains more exposed to fire.
- Fast-growing pine boosts yields but builds fuel; slower oak limits fire losses.

# 1. Introduction

Emergent wildfire risk challenges forest management planning, as climate-driven changes in fire extent, intensity, and frequency transform previously safe managed forests into vulnerable ones and threaten long-term ecosystem service provision (Girardin et al., 2013; Patacca et al., 2023; Tapias et al., 2004). Recently, severe wildfires have affected the northwestern parts of Europe, the Pacific Northwest of the United States, the boreal Canadian and Eurasian forests, and Southern Amazonia in Brazil (Cardíl et al., 2023; Davis et al., 2017; Fernandez-Anez et al., 2021; Senande-Rivera et al., 2022). In Europe, more than half a million (504,002) hectares were burnt in 2023, trailing 2017 (988,427 ha) and 2022 (837,212 ha). Across the continent, wildfires accounted for 24% of total timber damage between 1950 and 2020, equivalent to $12.5 \times 10^6$ $m^3$ $yr^{-1}$ (Patacca et al., 2023), with projections indicating this could rise over the next 15 years (Grünig et al., 2026; Lanet et al., 2024; Seidl et al., 2014). Beyond direct economic losses, fires are reshaping forest composition, altering soil properties, reducing biodiversity, and increasing vulnerability to further disturbances such as pest outbreaks and erosion (Pausas & Keeley, 2019).

Fire behavior emerges from the interaction of climate, vegetation, management practices, and human activities—making prediction inherently uncertain (Costafreda-Aumedes et al., 2017; Curt et al., 2016; Hantson et al., 2016; Puettmann, 2011). Vegetation type, fuel availability, and forest management practices further influence fire regimes, while human activities—such as land-use changes, fire suppression strategies, and ignition sources—add another layer of complexity. These complex interdependencies challenge efforts to accurately predict future wildfire trends and their long-term impacts on ecosystems and society.

To address this complexity, a wide array of fire models has been developed. High-resolution empirical and physical models (e.g., FARSITE, FIRETEC) provide detailed insights into fire behavior, but their computational demands limit scalability (Finney, 1998; Morvan et al., 2006; Silva et al., 2022). In contrast, modules integrated to dynamic global vegetation models (e.g., CLM-Li, INFERNO) offer broader-scale forecasts tied to climate and land cover (Kloster et al., 2010). Finally, recent process-based models like MC-FIRE, Reg-FIRM, CTEM, SPITFIRE, LMfire or LPX adapted fire forecasting to describe small-scale fire dynamics and scale up to the climate grid cell unit (Hantson et al., 2016; Pfeiffer, 2013). These models integrate components such as fire intensity, combustion completeness, and vegetation mortality, building on foundational work like the Rothermel fire spread model (Andrews, 2018; Burgan & Rothermel, 1984).

Meanwhile, process-based fire models have advanced our understanding of fire impacts on forest dynamics and ecosystem services, yet their potential to inform strategic management decisions remains largely untapped. A few optimization methods currently account for the uncertainty and variability introduced by such environmental extremes (Botequim et al., 2017; Couture & Reynaud, 2011; Finney et al., 2007; Garcia-Gonzalo et al., 2014; Spies et al., 2017; Taylor et al., 2013). These frameworks combined spatially explicit timber harvest with fire occurrence decision models, fire suppression and spread models, and fire protection value models (Acuna et al., 2010). Most existing methods focus primarily on immediate, stochastic losses of ecosystem services, failing to capture the full impact of fire behavior over time. This results in the under-representation of key long-term impacts—such as delayed mortality, shifts in regeneration, and cumulative ecosystem changes—which, in turn, limit strategic planning (Pasalodos-Tato et al., 2013; Qi & Zhuang, 2024).

In contrast, robust optimization addresses uncertainty by evaluating the gap between planned expectations and real-world outcomes. With minimal detailed input data, it assesses how landscape management compositions perform under various scenarios and implementation errors (e.g., prediction inaccuracies, productivity fluctuations) (Knoke et al., 2015, 2025; Reith et al., 2022). Combining robust optimization with process-based fire models thus offers new avenues to better support forest management planning under fire risk.

Pine plantation forests are vital for global wood production, yet highly vulnerable to wildfires, where management strategies (species selection, silvicultural regimes) influence fire vulnerability by altering fuel loads and moisture. We developed and coupled a process-based fire model (GO+Fire) with robust optimization to test whether this integrated approach can support reliable decision-making under fire uncertainty, and whether increasing the share of broadleaves in pine-dominated landscapes improves fire resilience. Specifically, we ask: (1) Is coupling robust optimization with a process-based model an effective method for managing fire risk? (2) Are broadleaf species less vulnerable to fire productivity loss compared to pine species? (3) Which ecosystem service objectives (timber, carbon) can current management practices effectively protect under fire risk?

# 2. Material and Methods

We selected the Landes de Gascogne Forest as a case study to assess fire risk in plantation systems and potential strategies to reduce future fire impacts. In plantation-dominated regions, wildfires can disrupt local economies, supply chains, and ecosystem services. Emergent fire activities (Figure S1) increased the need to balance productivity with risk reduction. In 2022, wildfires burned more than 60,000 ha in the Landes de Gascogne, Europe's largest contiguous plantation forest, highlighting the vulnerability of uniform pine stands. This event renewed interest in whether diversification through deciduous species or mixed stands could improve ecological and operational resilience (Carcaillet et al., 2022; Huffman et al., 2020).

We integrated a fire module into a process-based growth model (GO+) to simulate the fire response of Maritime Pine and Pedunculate Oak species in the LdG region. We then developed an optimization analysis to assess the impact of fire risk on forest management planning strategies.

## 2.1. Integration of a fire bio-geochemical module into GO+

We improved the forest growth and management model GO+ to simulate interactions among management regimes, species selection, fire behavior, and local climate. By integrating fire-spread equations (Peterson & Ryan, 1986; Pfeiffer, 2013; Rothermel, 1972; Wagner, 1973) into GO+ (Figure 1), we simulated both immediate fire impacts and their cascading effects on forest functioning and service provision.

GO+ simulates the growth of monospecific forests under different climate conditions and management regimes (Moreaux et al., 2020) on a daily basis. The management regimes include: (i) a reference regime for each species, (ii) short- and long-rotation regimes, (iii) a drought-adaptive regime in which thinning is triggered preventively to limit water stress, (iv) density management regimes based on a relative density index (RDI), which expresses stand density as a proportion of a species-specific maximum and is used to maintain stands within predefined “open” or “dense”, (v) coppice management for oak, and (vi) self-thinning regimes where stand density is controlled only by natural mortality (Table S1). These regimes mainly differ in regeneration type, stem density, thinning intensity and timing, and harvest schedule and product assortment.The Rothermel equations simulate fire events using forest structure, fuel type, dryness, and climate variables (Thonicke et al., 2010). These equations were previously integrated into the SPITFIRE model and coupled with large-scale vegetation models, including LPJ (Pfeiffer, 2013), JSBACH (Lasslop et al., 2014), and ORCHIDEE (Yue et al., 2014). We integrated them into GO+, refining fuel moisture, fire intensity, and damage evaluation to better capture fire behavior in maritime pine and sessile oak (Figure S2, Appendix A). The resulting GO+Fire model simulates daily changes in fuel types, fire susceptibility, fire spread, crown scorching, and plant mortality. It does not include ignition sources but accounts for weather-related fire extinction.

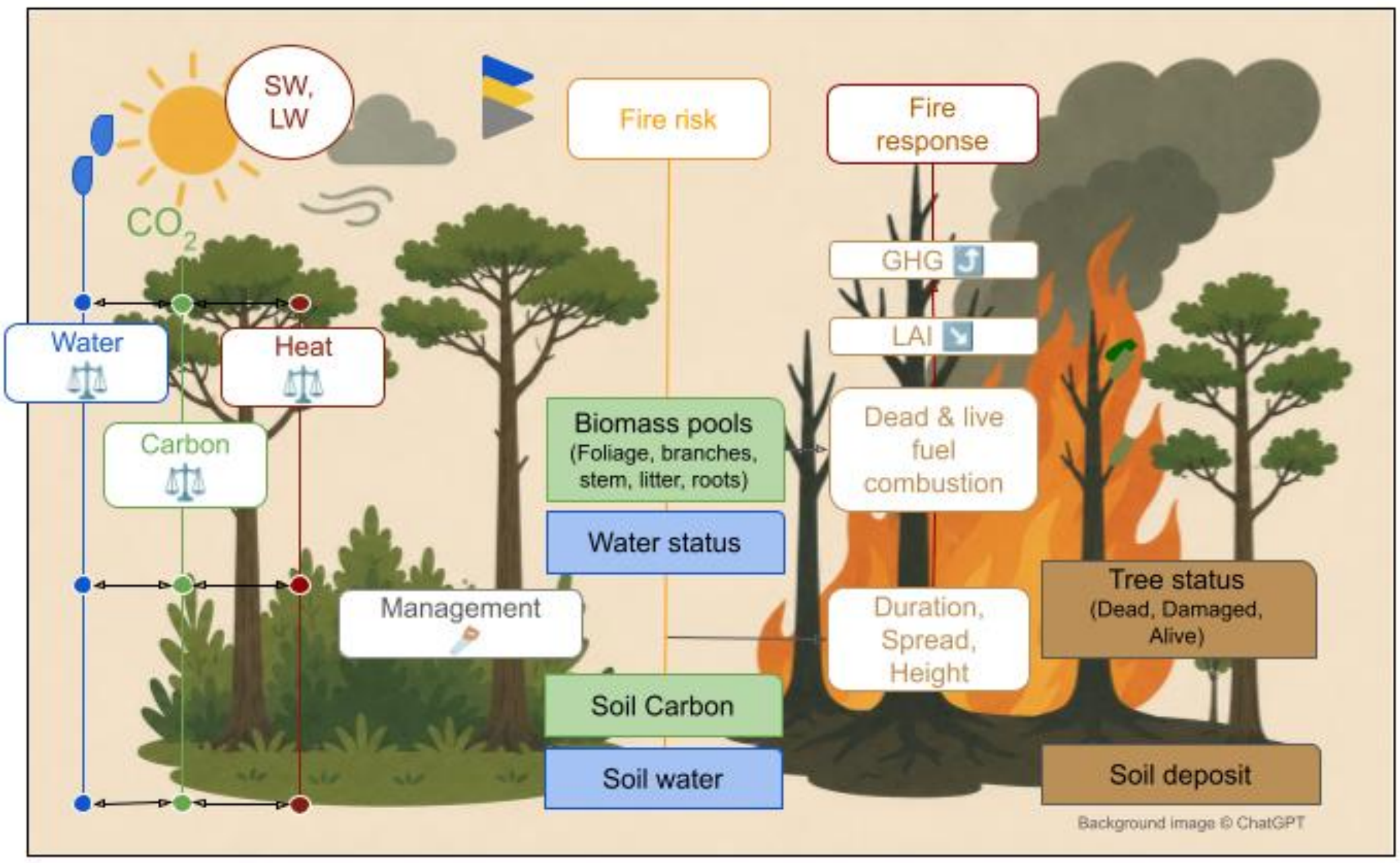

*Figure 1 Conceptual Overview of Fire Exposure and Response in GO+ Forest Model.* The diagram outlines the key processes linking climate inputs (radiation, heat fluxes, wind, precipitation) to ecosystem functioning (trees, understorey, soil), fire risk evaluation, and the subsequent damages from fire spread and intensity. Atmospheric fluxes (e.g., SW and LW – shortwave and longwave radiation; Heat – sensible and latent heat flux; $CO_2$ – carbon dioxide exchange) are affecting the phenology and growth of trees and vegetation. These processes are also regulated by the selected management regime. The biomass thus accumulated in the forest is distributed among stems, branches, foliage, coarse roots, and understorey, and is categorized into fuel classes (1h, 10h, 1000h, 10000h dead fuel, and live fuel). The water status of these fuel classes is calculated based on climate parameters. The quantity of fuel, its moisture content, and the climatic conditions determine the fire-risk exposure of the forest and the potential fire characteristics (duration, spread, flame scorch height). The effects of fire on vegetation compartments conditions mechanisms such as greenhouse gas emissions, soil carbon pool dynamics, tree mortality, decline (with leaf area index (LAI) decrease), and growth slowdown.

We evaluated GO+Fire by simulating 40 historical fires in the Landes de Gascogne under a standard 20-year-old pine stand and comparing hourly burned area with historic records (Figure S3) and wildfires >20 ha (Figure S4) from the BDIFF database (https://bdiff.agriculture.gouv.fr/) (IGN, 2025). GO+Fire underestimated mid-intensity fires but reproduced the observed range of extreme spread (≈2 to >1,000 ha $h^{-1}$), supporting its use for scenario exploration.

## 2.3. Simulations under and without fire risk

Simulations were conducted for two climate change scenarios (RCP 4.5 and RCP 8.5), three 30-year climatic periods (present, near future, and far future), two risk scenarios (no risk and fire risk), and two geographic sites of the LdG forest to contrast temperature and precipitation ranges (hydric constrained northern site, annual P = 750 mm/yr and unconstrained piedmont site, annual P = 1400 mm/yr). These directly affect the net primary productivity, production type, and the abundance, composition, and moisture content of biomass fuel (Hantson et al., 2016). Each simulation was initialized with an equally balanced age-class of pine forests (e.g., 25% of 2-, 16-, 30- and 43-year-old stands) and a soil available water content (AWC) distribution representative of the LdG region (40% of area at 40 mm, 40% at 80 mm, and 20% at 120 mm available water content; Arbez et al., 2017).

Fire-risk simulations included one fire per 30-year period. For each period, the ignition date was chosen within the middle decade of the period, e.g. 2016-2026 for the first period 2006-2036, and between the second and third quartiles of the longest period where the FDI value remained above 0.85 (Table S3). For example, the date chosen for the period 2006-2036 was the day 228 of year 2021 where the FDI remained above 0.85 from DOY 204 to DOY 240. The ignition date thus determined was applied consistently across all management regimes, sites, and climate scenarios. This enabled us to compare their performance under peak fire danger.

## 2.4. Provisions and fire damage indices

Five pine and seven oak management regimes, differing in thinning intensity, timing, and rotation length (Table S2), were simulated under two climate scenarios (RCP 4.5 and RCP 8.5) and three 30-year periods. Outcomes included in situ carbon stock and greenhouse gas (GHG) budget (g C or $CO_2$eq $m^{-2}$ over 30 years), harvested timber and stand productivity ($m^3$ $ha^{-1}$ $yr^{-1}$), and fire damage indices.

After building the simulation framework, the analysis focused on forest service provision and fire damage indices for two species—maritime pine and oak—at two sites with contrasting precipitation regimes. We analyzed the following variables:

(i) In situ carbon storage is the net change in ecosystem and soil carbon over each 30-year period.

(ii) The GHG budget combines net ecosystem exchange with fire-related emissions ($CO_2$eq, including $NO_x$, CO, $CH_4$, VOCs, and particulates) and carbon exported in harvested wood (potential delayed release of GHG). Positive GHG budget indicates net GHG uptake by the ecosystem, negative budget indicates net emissions.

(iii) Productivity refers to the biological production of wood and accounts for annual biomass growth and timber harvests (Moreaux et al., 2020).

(iv) The damage index quantifies the loss in provision caused by a fire ignited in each 30-year period. It is computed over the burned area and then multiplied by the probability of experiencing a similar event within that 30-year period (see Appendix 1 for details).

Regime performance under no fire is summarized in parallel-coordinate plots pooling both sites (Figure 2; Figure S5 which illustrate baseline variability in ecosystem services across RCP scenarios, climatic periods and species. For pine, timber and productivity range from 2.5 to 7 $m^3\ ha^{-1}\ yr^{-1}$, carbon storage is around 280,000–430,000 gC $m^{-2}$ over 30 years, and GHG fluxes range from 1 to 3 gC $m^{-2}\ yr^{-1}$. For oak, timber ranges from 0 to 6 $m^3\ ha^{-1}\ yr^{-1}$, productivity from 1 to 3.8 $m^3\ ha^{-1}\ yr^{-1}$, storage from 300,000 to 600,000 gC $m^{-2}$, and GHG from −8,000 to 0 gC $m^{-2}$. We further explained these patterns by running a multivariate random forest analysis (MultiOutputRegressor with RandomForestRegressor; Pedregosa et al., 2011) to quantify the effects of climate scenario, climatic period, management regime, species, available water content, site, and year on damage at harvest, GHG budget, carbon storage, productivity, and fire probability $f_{fire}$, defined in Appendix 1.

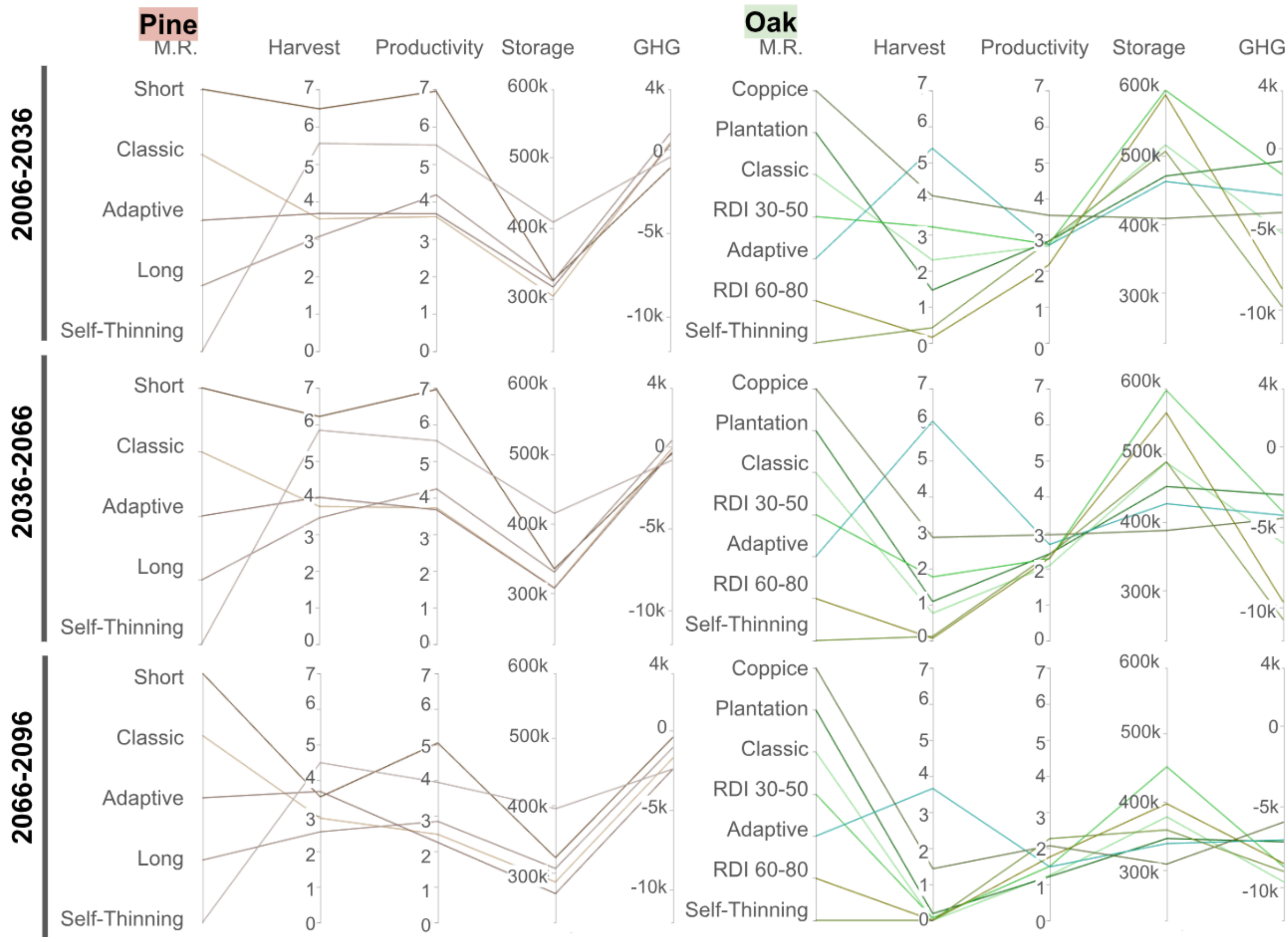

*Figure 2 Performance profile of forest management regimes under RCP 8.5. This parallel coordinate plot illustrates the performance of different management regimes across key indicators: Harvest ( timber production) and forest productivity ($m^3\ ha^{-1}\ yr^{-1}$), as well as carbon storage and ecosystem greenhouse gas budget( GHG buget in g $CO_2$eqe $m^{-2}$ over 30 years). Each management regime was simulated over two sites with contrasting rainfall regimes, and their climate responses were averaged. Results are presented for three climatic periods: 2006–2036 (A), 2036–2066 (B), and 2066–2096 (C), and for two species: Pinus pinaster and Quercus robur.*

## 2.5. Robust optimization

Process-based simulations under two climate-change scenarios assessed forest responses over 30-year periods and defined the range of possible fire-related outcomes, from no-fire conditions to highly severe events. Robust optimization then identified management portfolios that perform satisfactorily across this uncertainty space (Knoke et al., 2025). Rather than representing every fire scenario explicitly, robust optimization focuses on extreme outcomes to identify management regimes less sensitive to fire variability. It yielded portfolios that were more robust under both favorable and adverse conditions. The objective was to minimize the maximum relative regret $\beta$ associated with selecting one management landscape composition over another in terms of provision loss and damage risk. Robust optimization proceeded in two steps.

First, we defined the uncertainty space of each provision outcome and fire damage. For each management regime, we quantified the interannual variability of provision outcomes and the variability of fire damage to those provisions. Provision variability was estimated on a simulation without fire risk by aggregating annual outputs over AWC and age classes using their respective landscape weights, and then computing the standard deviation over the 30-year period. For timber and productivity, the standard deviation was calculated directly at the forest scale. For each provision p (carbon sequestration, carbon flux, timber production, or productivity) and management regime ***i***, $\sigma_1$ denotes the standard deviation of the 30-year series of provision outcomes around its expected value $E(c_{p,i})$ as follows:

*Equation 1*

$$\sigma_1 = \sqrt{\left(\frac{\sum_1^T (c_{p,i,t} - E(c_{p,i}))^2}{T}\right)}$$

with *T* corresponding to the 30 years of a climatic period (2006-2036, 2036-2066 and 2066-2096).

Since simulated fire events do not repeat over the 30-year period, we applied Taylor's law to estimate the standard deviation of provision loss under fire as follows (Senf et al., 2025):

*Equation 2*

$$\sigma_2 = \sqrt{0.23 \times (l_{p,i})^{2.28}}$$

Here, $l_{p,i}$ is the damage index (mean disturbance rate in Senf et al.) of provision *p* under management regime *i*, and $\sigma_2$ represents the expected standard variation of fire disturbance rate to that provision. This formulation was calibrated by Senf et al. (2025) using a European wildfire disturbance dataset.

Based on these metrics, we identified a best-case and worst-case outcome (fire damage) that represents extreme but plausible realizations of uncertainty for each management regime, as follows:

- For provision outcomes

*Equation 3 : Best case i.e. expected performance*

$$E(c_{\{p,i\}}) \quad \text{for all provision } p \text{ in } P, \text{management regime } i \text{ in } I$$

*Equation 4 Worst case i.e. minimum performance*

$$E(c_{p,i}) - m * \sigma_1$$

$$for\ all\ p\ in\ P, i\ in\ I$$

- For fire damage index

*Equation 5 Best case i.e. expected danger index*

$$\mathrm{L}_{p,i}$$
$$for\ all\ p\ in\ P, i\ in\ I$$

*Equation 6 Worst case i.e. maximum danger index*

$$\mathrm{L}_{p,i} + m * \sigma_2$$

$$for\ all\ p\ in\ P, i\ in\ I$$

With m being the size of the uncertainty space, which for the purposes of this methodological paper we fix to 1 for all provisions and optimization experiments.

We finally constructed a finite set of uncertainty scenarios *U* by combining these best and worst outcomes (fire damage) across all management alternatives, including scenarios where all alternatives are favorable, all are adverse, and mixed configurations. We denote by $k_{p,u,i}$ the realization of provision *p (*outcome or damage index) under management regime *i* in uncertainty scenario *u*.

Second, we used this set of uncertainty scenarios to identify robust management portfolios. We aimed to find mixtures of management alternatives, expressed as percentage allocations, that maintain good performance under changing conditions. Robustness was assessed by testing each portfolio across all uncertainty scenarios, evaluating all plausible future forest responses to management and climate change. Optimization was carried out for each period and climate scenario, comparing allocations with and without fire risk, and exploring mixed landscapes combining pure pine and oak stands.

## 2.6. Optimization with climate change variability

We first optimized management portfolios under inter-annual variability. To identify combinations of management regimes that jointly maximize expected provision and minimize variability, we reformulated the optimization problem. We replaced the classic objective function by a system of robustness constraints that enables comparison and evaluation of portfolio performance across the full set of uncertainty scenarios (Knoke et al., 2025).

### Constraints:

The decision variable $X_{p,u,i}$ determines which area of the site is to be allocated to management alternative *i* under provision *p* and scenario *u*. The selection process followed these principles: Each provision *p* has multiple possible outcomes *k* across different management alternatives *i*. For a given scenario *u*, the selected outcome should be the closest to the maximal achievable outcome under that scenario comparing all possible management regimes. The relative performance of a provision *p* under scenario *u* is given by:

*Equation 7*

$$Q_{p,u} = \frac{\sum_i^I X_{p,u,i} * k_{p,u,i} - \min_i k_{p,u,i}}{\max_i k_{p,u,i} - \min_i k_{p,u,i}}$$

This formula ensured that performance values are normalized between 0 and 1, allowing fair comparison across different provisions. The numerator represents the selected outcome's deviation from the worst-case performance. The denominator represents the total possible performance range under scenario $u$.

The worst relative performance across all provisions should be minimized:

*Equation 8*

$$100 - Q_{p,u} \leq \beta$$

This ensured that no provision falls too far behind its best possible outcome, promoting balanced performance across all forest goods.

### Objective

The optimization problem sought to minimize the worst-case relative deviation across all provisions, i.e. the maximum relative regret:

*Equation 9*

$$Minimize\ \beta$$

## 2.6. Optimization with climate change variability and fire risk variability

We applied equations 2,5-9 with k representing fire damage indices to consider the variability of fire damages.

We ran optimizations both for each provision individually and for all provisions jointly. For each case, we compared solutions with fire risk (Section 2.6) and without fire risk (Section 2.7). Robust optimizations were conducted in a multi-species setting, where the decision space included both species choice (Maritime pine vs. Sessile–pedunculate oak) and their associated management regimes, yielding 12 decision alternatives in total (Table S1). All models were solved using CPLEX 22.1.1.0 (IBM Corp., 2019) via linear programming.

# 3. Results

## 3.1. Influence of fire on the vulnerability of tree species

Understanding species-specific fire responses is essential for assessing how fire may affect forest services under climate change. Figure 3 shows the difference between the no-fire-risk and fire-risk scenarios (no fire risk minus fire risk). Positive values, therefore, indicate higher provision under no fire risk, whereas negative values indicate greater provision under fire risk, except for the GHG budget, where positive values indicate lower emissions under fire risk. Significantly larger losses in timber production, carbon storage, and the GHG budget were found for pine than for oak, indicating higher vulnerability. On average, pine lost about 1 $m^3$ $ha^{-1}$ $yr^{-1}$ of timber production under fire risk, whereas oak showed negligible losses, although losses reached up to 4 $m^3$ $ha^{-1}$ $yr^{-1}$ under specific conditions. In worst-case scenarios, pine experienced the greatest risks for timber production and GHG emissions, while oak showed the largest potential losses in carbon storage.

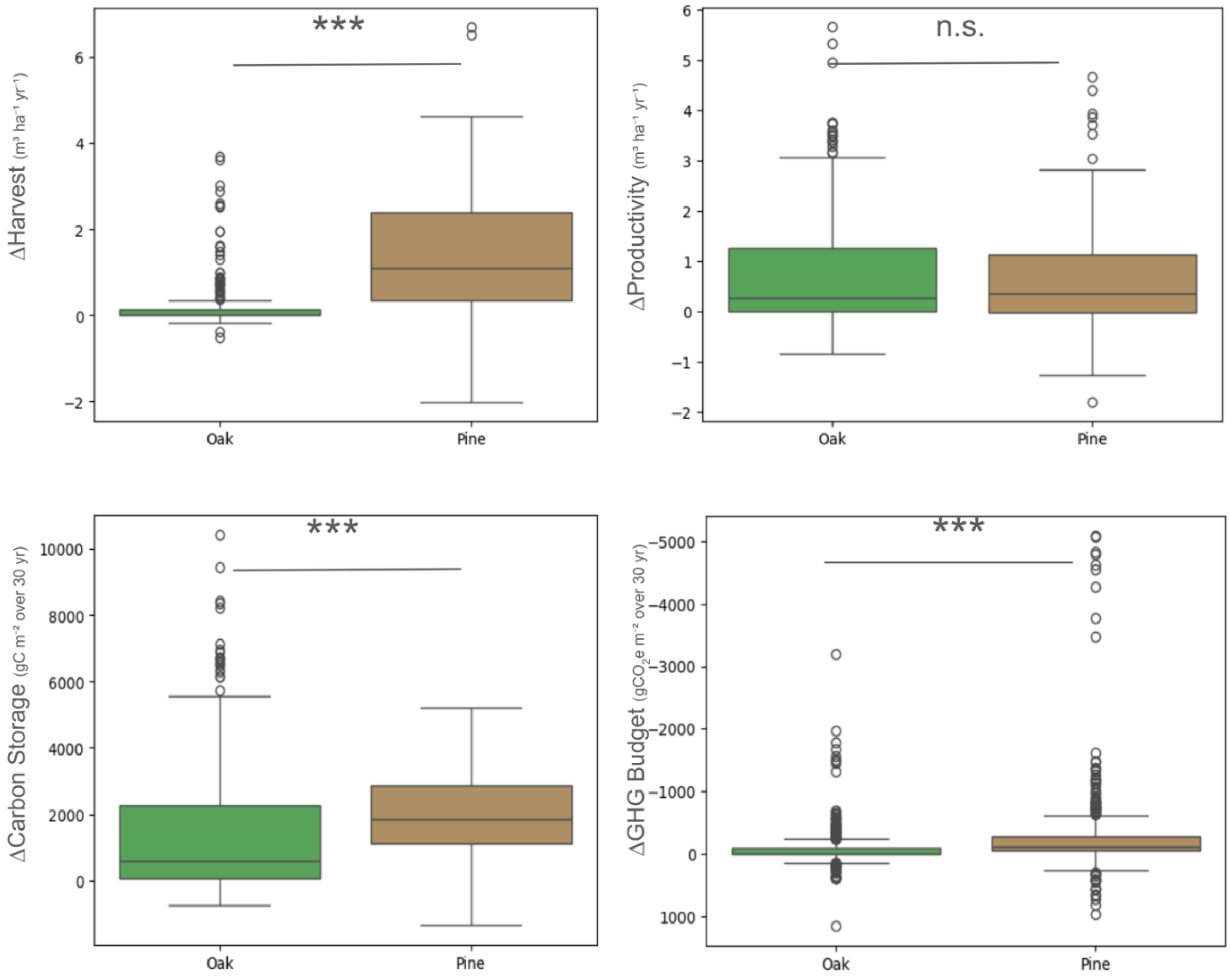

*Figure 3 Projected distribution of fire damages for four ecosystem services for Oak and Pine management regimes.* We assess variation in fire impacts on timber provision, productivity, carbon storage, and greenhouse gas exchange for oak and pine species. Each management regime is stress-tested under no-fire and fire-risk scenarios. The figure shows the magnitude of change, defined as the provision outcome under no fire risk minus the provision outcome under fire risk. Five management regimes were evaluated for pine and seven for oak. These performance differences are analyzed across two RCPs (4.5 and 8.5), three 30-year climate periods (2006–2036, 2036–2066, 2066–2096), and two site conditions (hydric-constrained and unconstrained). The resulting distributions combine all combinations of RCP, site, period, and management regimes to particularly highlight variability across species. Above each plot, the symbols provide the statistical significance of the difference between oak and pine (*** $p < 0.001$ or n.s. $p = 0.684$).

### 3.2. Management regimes drive the difference in fire-response variability

Multivariate random forest analysis identified management regime as the main driver of fire damage variability, followed by site and year (importance scores: 0.36, 0.32, and 0.14; Table S4). Figure 4 illustrates how fire damage varied among management regimes across climate scenarios, periods, and sites for each ecosystem service provision.

Figure 5 ranks management regimes by their average loss of provision due to fire. Carbon storage exhibited the widest range of average management performance (13000 to 38000 g C $m^{-2}$ over 30 years), whereas mean harvest losses varied comparatively little among management regimes, although short-rotation management presented a wide distribution. The ranking of management regimes differed among ecosystem services (Figure 4). In pine, self-thinning produced the greatest timber and carbon losses, whereas short-rotation management showed the highest losses in productivity and GHG balance, likely reflecting the greater fire sensitivity of younger stands. In oak, coppice management generated high losses in carbon storage and productivity but comparatively low losses in GHG balance. These contrasting rankings show that no management regime consistently minimizes fire damage across all ecosystem services and scenarios.

High-yield management regimes frequently experienced large fire losses, although this pattern was not universal (Figure 4). Pine short-rotation management sustained the highest productivity losses (0.5 $m^3$ $ha^{-1}$ $yr^{-1}$), offsetting part of its high production potential (6.7 $m^3$ $ha^{-1}$ $yr^{-1}$). Likewise, oak drought-adaptive management showed the greatest timber vulnerability despite being among the most productive oak strategies. In contrast, pine self-thinning exhibited the largest carbon losses (up to 6000 g C $m^{-2}$ over 30 years) despite relatively low carbon stocks. Overall, fire vulnerability did not follow a simple intensive-to-extensive management gradient but instead depended on the combination of production potential and stand structural characteristics.

Forest management regimes inherently differed in their capacity to provide ecosystem services under undisturbed conditions (shaded bars, Figure 4). When fire was introduced, these initial differences were amplified, increasing both the magnitude and variability of ecosystem service losses. The response depended on the ecosystem service considered, reflecting differences in stand structure and management characteristics. For example, carbon storage was highly sensitive to fire under the oak adaptive regime, whereas timber production remained comparatively stable (Figure 4). These findings highlight that management regimes differ in their effectiveness under fire risk and that their performance also depends on the combination of climate scenario, climatic period, and site conditions. This context dependence creates multiple response pathways that challenge robust forest planning.

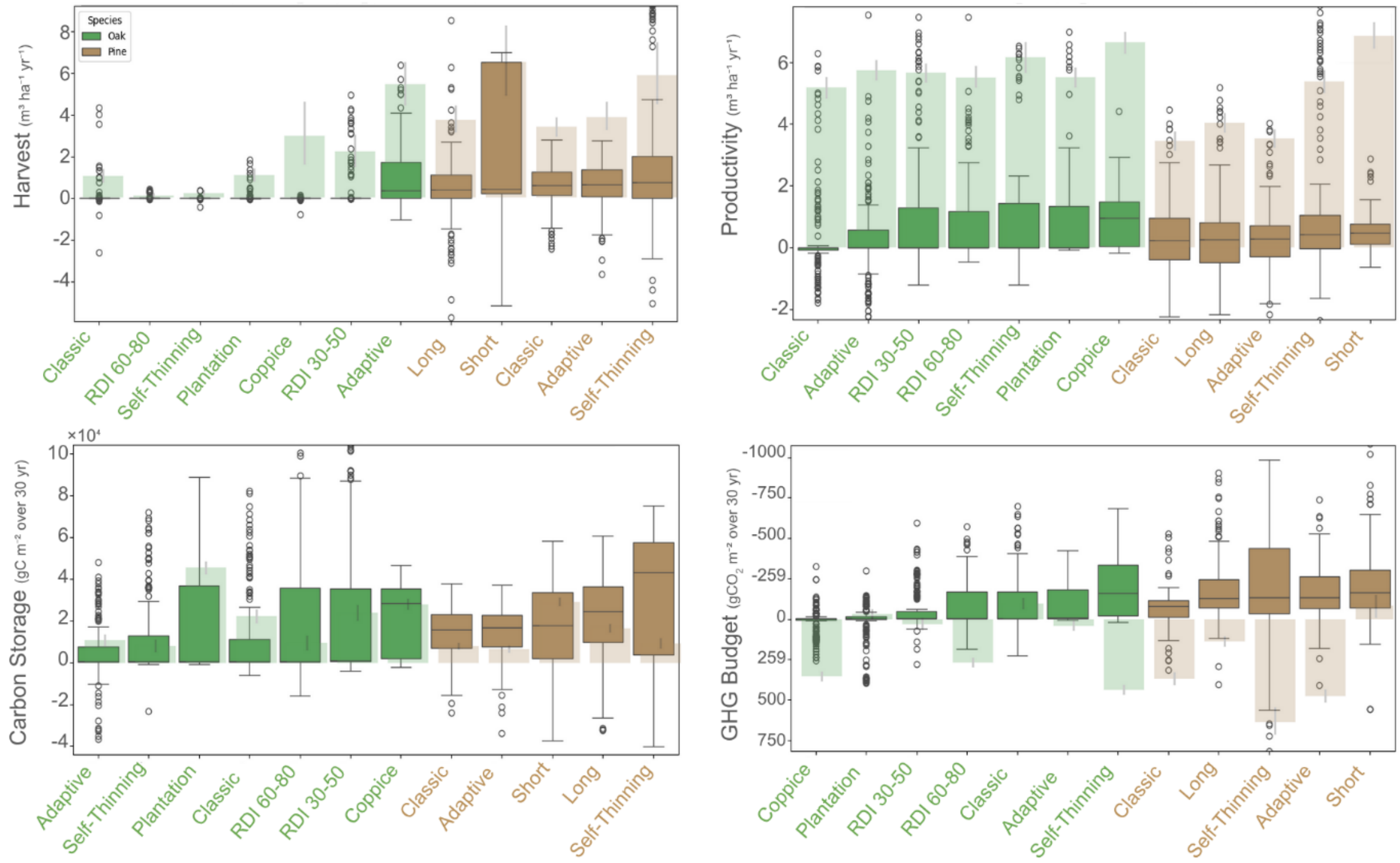

*Figure 4 Fire damage to the production of four ecosystem services under different management regimes.*

We compare baseline ecosystem service provision (shaded bars) and vulnerability to fire (plain box plots) across management regimes. Shaded bars show the mean provision under no-fire scenarios, averaged across two RCPs (4.5 and 8.5), three climate periods (2006–2036, 2036–2066, and 2066–2096), two study sites, and their corresponding soil types. Error bars indicate the standard error of the mean. Plain box plots show the distribution of fire damage, calculated as the difference in provision between paired fire and no-fire scenarios under the same climatic and site conditions. They therefore capture variation in fire damage across RCPs, climate periods, sites, and soil types. For each ecosystem service and species, management regimes are ranked by increasing median fire damage. Figure-reading keys are provided in Figure S6.

## 3.3. Fire affects forest management planning strategies

We compared landscape portfolios optimized by the robust method with and without fire. Figure 5 showed the optimization results for the site constrained by hydric conditions, under the RCP 8.5 projection for the far future (2066–2096). Under these conditions, this site exhibited the greatest fire risk. In the absence of fire risk, pine self-thinning emerged as the most productive option for timber production, alongside pine drought-adaptive management, an intermediate choice. Oak self-thinning was the best strategy for promoting productivity, outperforming oak coppice and pine short-rotation management, both of which were also considered viable candidates. Additionally, oak self-thinning supported a better greenhouse gas (GHG) budget than pine self-thinning and adaptive management strategies. These results were consistent with the average production levels observed in Figure 4. The preference for intermediate strategies can be attributed to their more stable performance over the simulated years.

When fire risk was incorporated, single-objective optimization revealed that management regime selection was influenced by trade-offs between performance and fire vulnerability (Figure 5). Under the timber provision, two pine adaptive and classic management regimes were selected. Despite their high fire-damage potential, they offered high timber yield. Conversely, the oak relative density index (RDI) 60 -80 balanced this risk by presenting no fire damage but limited timber provision. Despite these trade-offs, recommendations for productivity remained largely unchanged from the no-fire planning scenario. No significant differences in fire-risk response were found between oak and pine species, supporting the continued selection of the most productive options. For the GHG budget, oak adaptive management was again included to mitigate high emission risks. When carbon storage became the priority, a diversified portfolio of management strategies emerged as the most effective approach Figure 5.

Multi-objective optimization under fire risk similarly showed a pattern of diversification; however, overall carbon storage performance remained limited, securing only 20% of the target under fire-risk conditions. This highlighted the challenge of fully securing ecosystem services, particularly carbon storage, in fire-prone landscapes, even with a diversified and optimized management portfolio.

These findings underscore the cost of introducing fire risk and its inherent uncertainty into management planning. Robust optimization selects management regimes that ensure satisfactory ecosystem service provision, but this often comes at the expense of missed opportunities for higher returns. In other words, to secure provisions from potential fire events, decision-makers must accept deviations from the management regime that would otherwise yield the best outcomes. For timber production and GHG budget, fire-robust planning could guarantee 60% of the performance range between the minimum and maximum provisions across uncertainty scenarios u. Highly protective management options with lower returns were included, such as oak RDI 60–80 or adaptive management regimes. Interestingly, for carbon storage, the robust solution only guaranteed 20% of the performance range, suggesting that no combination of management regimes can fully immunize against fire risk.

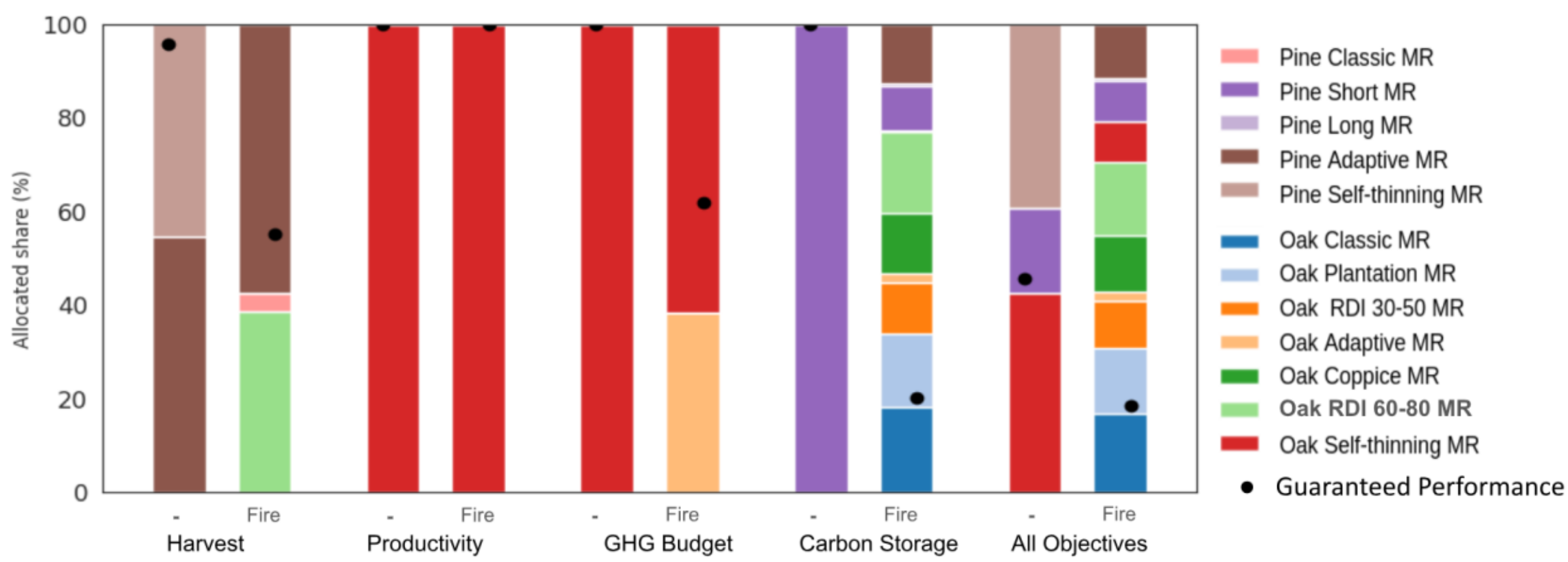

*Figure 5 Optimized allocation of management regimes to minimize regret under no-fire and fire scenarios for the period 2070-2100 under RCP 8.5.* Maximal regret is defined as the difference between the maximum and minimum provision performance across uncertainty scenarios. Robust optimization minimizes this regret by selecting management portfolios that reduce the deviation from the maximum attainable performance. The resulting **guaranteed performance** (black dots) represents the lowest provision performance maintained across all uncertainty scenarios for the selected portfolio. It is normalized between the minimum and maximum attainable performance, with 0% corresponding to the minimum and 100% to the maximum. Thus, a guaranteed performance of 20%, for example, indicates that the selected portfolio secures a performance level 20% of the way between the minimum and maximum attainable values. Provision objectives are harvest, productivity, GHG budget, and carbon storage, optimized under no-risk (–) and fire-risk (fire) scenarios.

## 4. Discussion

Effective forest management planning is increasingly challenged by rising fire risk and uncertainty in ecosystem service provision under climate change. Traditional approaches often fail to capture the nonlinear and interacting effects of fire, vegetation dynamics, and management interventions. Here, we developed a coupled modeling framework integrating the GO+ vegetation model with Rothermel fire behavior equations to evaluate the robustness of management–species strategies under uncertain fire conditions. By stress-testing alternative regimes across productivity, carbon storage, greenhouse gas emissions, and timber production, our framework identifies vulnerabilities and opportunities for risk-informed forest planning. We show that carbon storage is particularly sensitive to fire risk, driving the need for diversified and innovative management strategies (Özkan & SERENGİL, 2025; Santana et al., 2016). Extending previous diversification approaches (Knoke et al., 2016), we demonstrate that resilience depends not only on increasing diversity but also on ensuring that management options differ in their ability to reduce risk. This framework reveals how disturbances reshape optimal decisions and helps identify ecosystem services most vulnerable to underperformance (Knoke et al., 2025; Puettmann, 2021).

Our key contribution lies in integrating fire disturbance into forest planning through robust optimization. This approach accounts for uncertainties that cannot be fully captured by stochastic processes (Gabrel et al., 2014; Labarre et al., 2025; Reith et al., 2022), incorporating variability in fire occurrence and damage without simulating individual ignition events or fire sequences. By representing fire as a probabilistic disturbance, the model captures variation in expected damage across climate scenarios, periods, sites, and stand ages. Despite this simplified representation, fire inclusion substantially alters management decisions, shifting strategies from risk-neutral to risk-averse planning. Optimal solutions increasingly combine productive but vulnerable regimes with strategies that secure ecosystem service provision. These shifts reflect the variability in fire damage simulated by GO+-Fire across species and management regimes (Figure 4), demonstrating that the magnitude of damage and uncertainty are key predictors of strategy performance (Fernández-Guisuraga et al., 2023). This highlights the value of integrating even simplified fire models into forest management optimization.

Beyond changes in strategy selection, fire damage significantly altered the distribution of species–management combinations, notably increasing the prevalence of less vulnerable species. This change reflects species-specific differences in fire response. The model captures these through detailed simulations of fuel composition and structural traits. For instance, *Pinus pinaster* is more prone to biomass loss due to higher flammability, yet benefits from fire-resistant traits such as thick bark and exfoliating structure, particularly in mature stands (Fernandes et al., 2008; Singh et al., 2022; Varner et al., 2022). In contrast, *Quercus robur* and *Q. petraea* exhibit lower immediate fire damage but suffer physiological impairments—such as reduced leaf area and xylem degradation—after repeated fire and drought exposure, which constrain long-term productivity (Botequim et al., 2017; Gričar et al., 2020; Lopes et al., 2024). These contrasting traits also shape service provision: fast-growing, frequently harvested pine accumulates more biomass and fuel, increasing potential fire losses, whereas slower-growing, less harvested oak maintains lower biomass stocks and thus lower loss rates (Cabrera et al., 2025).

Management regimes strongly influenced fire responses, particularly carbon storage, reflecting the sensitivity of understorey and litter biomass to fire disturbance. Because fire occurred mid-simulation, stands with greater pre-fire biomass experienced larger carbon losses. Following fire, the regeneration of young tree stands was instead faster than older ones. Carbon responses therefore depended on fire damage intensity, recovery dynamics, and thinning practices. In contrast, timber volume showed lower variability because it reflects discrete management interventions rather than total stand biomass (Figure 2). These differences should be interpreted considering that GO+-Fire evaluates regimes at stabilized states and even aged class distribution and does not account for the transition from current forest conditions. Consequently, part of the variation reflects differences in initial stand structure and age distribution rather than management design alone. For example, self-thinning regimes are initiated in older stands and with longer rotations, resulting in higher biomass stocks but lower carbon sequestration rates (Table S1). This equilibrium-based approach was selected to isolate the long-term stability and vulnerability of management regimes under climate change, independently of transition pathways.

With this in mind, our simulations reveal that management regimes play a critical role in modulating fire responses. Fernández- Guisuraga (Fernández-Guisuraga et al., 2022) and Fernandes & Rigolot (Fernandes & Rigolot, 2007) recommend managing surface fuel loads and structural continuity to reduce fire risk—through pruning, thinning, or ladder fuel removal—while noting that thinning alone can increase wind penetration and fuel drying. In line with these insights, our results show that rotation length and thinning intensity define fuel profiles and stand structure, ultimately influencing fire behavior. Notably, high-yield regimes often involve greater fire-exposure trade-offs. Pine short-rotation regimes, while highly productive, suffered the greatest fire-induced productivity losses (0.5 $m^3$ $ha^{-1}$ $yr^{-1}$). Similarly, oak adaptive management, despite stabilizing productivity, experienced the highest timber volume losses among oak regimes. These outcomes reflect not just structural design more prone to lethal crown fires and but also the combined effect of pre-fire biomass accumulation and post-fire management dynamics.

While this study focuses on fire as a dominant disturbance, it is essential to recognize that forest fires occur within a broader context of cascading and compound risks (Bastit et al., 2023). Climate change exacerbates forest vulnerability not only to fire, but also to drought and pest outbreaks, which interact in complex ways to shape ecosystem responses and management outcomes (Zscheischler et al., 2018). For instance, fire regimes may shift from fuel-driven to drought-driven, and post-fire bark beetle outbreaks can increase the severity and frequency of subsequent fire events (Fernández & Costas, 1999; Fernández-Guisuraga et al., 2023). The developed GO+-Fire model advances multi-risk assessments by capturing the influence of key climate variables on drought stress, fuel moisture, and forest productivity. Through this integration, the model enables a more holistic evaluation of forest resilience and management strategies. Such a comprehensive perspective is crucial to designing flexible and robust adaptation pathways in the face of accelerating environmental change (Puettmann, 2021).

In conclusion, explicitly accounting for fire risk substantially changes forest management decisions and the selection of optimal strategies. Robust optimization identifies management regimes that account for uncertainty while highlighting how management outcomes depend on species traits and fire-damage variability. These results provide a framework for more resilient, service-oriented forest planning in fire-prone landscapes. Future work should integrate robust optimization with dynamic process-based models (Knoke et al., 2020; Pérez-Romero et al., 2025) to capture emergent risk, management transitions, and long-term changes in forest structure and species composition, thereby producing recommendations that better reflect real-world forest dynamics.

# Reference

Acuna, M. A., Palma, C. D., Cui, W., Martell, D. L., & Weintraub, A. (2010). Integrated spatial fire and forest management planning. *Canadian Journal of Forest Research*, *40*(12), 2370‑ 2383.

Andrews, P. L. (2018). *The Rothermel surface fire spread model and associated developments : A comprehensive explanation*.

Bastit, F., Brunette, M., & Montagné-Huck, C. (2023). Pests, wind and fire : A multi-hazard risk review for natural disturbances in forests. *Ecological Economics*, *205*, 107702.

Botequim, B., Arias-Rodil, M., Garcia-Gonzalo, J., Silva, A., Marques, S., Borges, J. G., Oliveira, M. M., & Tomé, M. (2017). Modeling Post-Fire Mortality in Pure and Mixed Forest Stands in Portugal—A Forest Planning-Oriented Model. *Sustainability*, *9*(3), Article 3. https://doi.org/10.3390/su9030390

Burgan, R. E., & Rothermel, R. C. (1984). *Behave : Fire behavior prediction and fuel modeling system, fuel subsystem* (Vol. 167). US Department of Agriculture, Forest Service, Intermountain Forest and Range Experiment Station.

Cabrera, S., Alexander, H., Willis, J. L., Milton, T., & Anderson, C. (2025). *Fuel Loads and Composition Across a Pine Dominance Gradient within Unmanaged, Fire-Dependent Pine-Oak Mixedwoods of the Southeastern U.S* (SSRN Scholarly Paper N° 5143946). Social Science Research Network. https://doi.org/10.2139/ssrn.5143946

Carcaillet, C., Michalet, R., Fréjaville, T., & others. (2022). Les «Forêts» de pins maritimes d'Aquitaine, des nids à incendie. *The Conversation*, *30*.

Cardíl, A., Tapia, V. M., Monedero, S., Quiñones, T., Little, K., Stoof, C. R., & de-Miguel, S. (2023). Characterizing the rate of spread of large wildfires in emerging fire environments of northwestern Europe using Visible Infrared Imaging Radiometer Suite active fire data. *Natural hazards and earth system sciences*, *23*(1), 361‑ 373.

Costafreda-Aumedes, S., Comas, C., & Vega-Garcia, C. (2017). Human-caused fire occurrence modelling in perspective : A review. *International Journal of Wildland Fire*, *26*(12), 983‑ 998.

Couture, S., & Reynaud, A. (2011). Forest management under fire risk when forest carbon sequestration has value. *Ecological Economics, Special Section - Earth System Governance: Accountability and Legitimacy*, *70*(11), 2002‑ 2011. https://doi.org/10.1016/j.ecolecon.2011.05.016

Curt, T., Fréjaville, T., & Lahaye, S. (2016). Modelling the spatial patterns of ignition causes and fire regime features in southern France : Implications for fire prevention policy. *International Journal of Wildland Fire*, *25*(7), 785‑ 796.

Davis, R., Yang, Z., Yost, A., Belongie, C., & Cohen, W. (2017). The normal fire environment—Modeling environmental suitability for large forest wildfires using past, present, and future climate normals. *Forest Ecology and Management*, *390*, 173‑ 186.

Fernandes, P. M., & Rigolot, E. (2007). The fire ecology and management of maritime pine (Pinus pinaster Ait.). *Forest Ecology and Management*, *241*(1), 1‑ 13. https://doi.org/10.1016/j.foreco.2007.01.010

Fernandes, P. M., Vega, J. A., Jiménez, E., & Rigolot, E. (2008). Fire resistance of European pines. *Forest Ecology and Management*, *256*(3), 246‑ 255. https://doi.org/10.1016/j.foreco.2008.04.032

Fernández, M. M. F., & Costas, J. M. S. (1999). *Susceptibility of fire-damaged pine trees (Pinus pinaster and Pinus nigra) to attacks by Ips sexdentatus and Tomicus piniperda (Coleoptera : Scolytidae).*

Fernandez-Anez, N., Krasovskiy, A., Müller, M., Vacik, H., Baetens, J., Hukić, E., & Cerda, A. (2021). Current wildland fire patterns and challenges in Europe : A synthesis of national perspectives. *Air, Soil and Water Research*, *14*, 11786221211028185.

Fernández-Guisuraga, J. M., Marcos, E., & Calvo, L. (2023). The footprint of large wildfires on the multifunctionality of fire-prone pine ecosystems is driven by the interaction of fire regime attributes. *Fire Ecology*, *19*(1), 32.

Fernández-Guisuraga, J. M., Suárez-Seoane, S., Fernandes, P. M., Fernández-García, V., Fernández-Manso, A., Quintano, C., & Calvo, L. (2022). Pre-fire aboveground biomass, estimated from LiDAR, spectral and field inventory data, as a major driver of burn severity in maritime pine (Pinus pinaster) ecosystems. *Forest Ecosystems*, *9*, 100022. https://doi.org/10.1016/j.fecs.2022.100022

Finney, M. A. (1998). *FARSITE, Fire Area Simulator–model development and evaluation* (N° 4). US Department of Agriculture, Forest Service, Rocky Mountain Research Station.

Finney, M. A., Seli, R. C., McHugh, C. W., Ager, A. A., Bahro, B., & Agee, J. K. (2007). Simulation of long-term landscape-level fuel treatment effects on large wildfires. *International Journal of Wildland Fire*, *16*(6), 712‑ 727.

Gabrel, V., Murat, C., & Thiele, A. (2014). Recent advances in robust optimization : An overview. *European Journal of Operational Research*, *235*(3), 471‑ 483. https://doi.org/10.1016/j.ejor.2013.09.036

Garcia-Gonzalo, J., Pukkala, T., & Borges, J. G. (2014). Integrating fire risk in stand management scheduling. An application to Maritime pine stands in Portugal. *Annals of Operations Research*, *219*(1), 379‑ 395. https://doi.org/10.1007/s10479-011-0908-1

Girardin, M. P., Ali, A. A., Carcaillet, C., Gauthier, S., Hély, C., Le Goff, H., & Bergeron, Y. (2013). Fire in managed forests of eastern Canada : Risks and options. *Forest Ecology and Management*, *294*, 238‑ 249.

Gričar, J., Hafner, P., Lavrič, M., Ferlan, M., Ogrinc, N., Krajnc, B., Eler, K., & Vodnik, D. (2020). Post-fire effects on development of leaves and secondary vascular tissues in Quercus pubescens. *Tree Physiology*, *40*(6), 796‑ 809.

Grünig, M., Rammer, W., Senf, C., Albrich, K., André, F., Augustynczik, A. L. D., Baumann, M., Bohn, F. J., Bouwman, M., Bugmann, H., Collalti, A., Cristal, I., Dalmonech, D., De Coligny, F., Dobor, L., Dollinger, C., Espelta, J. M., Forrester, D. I., Garcia-Gonzalo, J., … Seidl, R. (2026). Climate change will increase forest disturbances in Europe throughout the 21st century. *Science*, *391*(6789), eadx6329. https://doi.org/10.1126/science.adx6329

Hantson, S., Arneth, A., Harrison, S. P., Kelley, D. I., Prentice, I. C., Rabin, S. S., & Yue, C. (2016). The status and challenge of global fire modelling. *Biogeosciences*, *13*(11), 3359‑ 3375.

Huffman, D. W., Floyd, M. L., Hanna, D. P., Crouse, J. E., Fulé, P. Z., Meador, A. J. S., & Springer, J. D. (2020). Fire regimes and structural changes in oak-pine forests of the Mogollon Highlands ecoregion : Implications for ecological restoration. *Forest Ecology and Management*, *465*, 118087.

IGN. (2025). *Base de Données sur les Incendies de Forêts en France (BDIFF)*.

Kloster, S., Mahowald, N. M., Randerson, J. T., Thornton, P. E., Hoffman, F. M., Levis, S., Lawrence, P. J., Feddema, J. J., Oleson, K. W., & Lawrence, D. M. (2010). Fire dynamics during the 20th century simulated by the Community Land Model. *Biogeosciences*, *7*(6), 1877‑ 1902.

Knoke, T., Biber, P., Schula, T., Fibich, J., & Gang, B. (2025). Minimising the relative regret of future forest landscape compositions : The role of close-to-nature stand types. *Forest Policy and Economics*, *171*, 103410.

Knoke, T., Paul, C., Rammig, A., Gosling, E., Hildebrandt, P., Haertl, F., Peters, T., Richter, M., Diertl, K.-H., Maria Castro, L., Calvas, B., Ochoa, S., Anabelle Valle-Carrion, L., Hamer, U., Tischer, A., Potthast, K., Windhorst, D., Homeier, J., Wilcke, W., … Bendix, J. (2020). Accounting for multiple ecosystem services in a simulation of land-use decisions : Does it reduce tropical deforestation? In *GLOBAL CHANGE BIOLOGY* (Vol. 26, Numéro 4, p. 2403‑ 2420). WILEY. https://doi.org/10.1111/gcb.15003

Labarre, C., Domec, J.-C., Andrés-Domenech, P., Bödeker, K., Bingham, L., & Loustau, D. (2025). Improving forest decision-making through complex system representation : A viability theory perspective. *Forest Policy and Economics*, *170*, 103384.

Lanet, M., Li, L., Ehret, A., Turquety, S., & Le Treut, H. (2024). Attribution of summer 2022 extreme wildfire season in Southwest France to anthropogenic climate change. *Npj Climate and Atmospheric Science*, *7*(1), 267. https://doi.org/10.1038/s41612-024-00821-z

Lopes, L. F., Dias, F. S., Fernandes, P. M., & Acácio, V. (2024). A remote sensing assessment of oak forest recovery after postfire restoration. *European Journal of Forest Research*, *143*(3), 1001‑ 1014.

Moreaux, V., Martel, S., Bosc, A., Picart, D., Achat, D., Moisy, C., & Loustau, D. (2020). Energy, water and carbon exchanges in managed forest ecosystems : Description, sensitivity analysis and evaluation of the INRAE GO+ model, version 3.0. *Geoscientific Model Development*, *13*(12), 5973‑ 6009.

Morvan, D., Dupuy, J. L., Rigolot, E., & Valette, J. (2006). *FIRESTAR: a physically based model to study wildfire behaviour*.

Pasalodos-Tato, M., Mäkinen, A., Garcia-Gonzalo, J., Borges, J. G., Lämås, T., & Eriksson, L.-O. (2013). Assessing uncertainty and risk in forest planning and decision support systems : Review of classical methods and introduction of new approaches. *Forest Systems*, *22*(2), 282‑ 303.

Patacca, M., Lindner, M., Nabuurs, G.-J., & Schelhaas, M.-J. (2023). *Significant increase in forest disturbances since 1950s.*

Pausas, J. G., & Keeley, J. E. (2019). Wildfires as an ecosystem service. *Frontiers in Ecology and the Environment*, *17*(5), 289‑ 295. https://doi.org/10.1002/fee.2044

Pérez-Romero, J., González-Sanchis, M., Blanco-Cano, L., & Del Campo, A. D. (2025). Development of a multi-objective decision support system for eco-hydrological forest management that quantifies and optimizes ecosystem services related to Carbon, Water, Fire-risk and Eco-resilience (CAFE). *Journal of Environmental Management*, *380*, 125103. https://doi.org/10.1016/j.jenvman.2025.125103

Peterson, D. L., & Ryan, K. C. (1986). Modeling postfire conifer mortality for long-range planning. *Environmental management*, *10*(6), 797‑ 808.

Pfeiffer, M. (2013). *Modeling Terrestrial Paleobiogeochemistry : Linking Fire, Biochemistry, Humans, and Climate* (Numéro 5721). EPFL.

Puettmann, K. J. (2011). Silvicultural Challenges and Options in the Context of Global Change : “Simple” Fixes and Opportunities for New Management Approaches. In *JOURNAL OF FORESTRY* (Vol. 109, Numéro 6, p. 321‑ 331). OXFORD UNIV PRESS INC.

Puettmann, K. J. (2021). Extreme Events : Managing Forests When Expecting the Unexpected. *Journal of Forestry*, *119*(4), 422‑ 431. https://doi.org/10.1093/jofore/fvab014

Qi, J., & Zhuang, J. (2024). A review of optimization and decision models of prescribed burning for wildfire management. *Risk Analysis*.

Reith, E., Gosling, E., Knoke, T., & Paul, C. (2022). Exploring trade-offs in agro-ecological landscapes : Using a multi-objective land-use allocation model to support agroforestry research. *Basic and Applied Ecology*, *64*, 103‑ 119. https://doi.org/10.1016/j.baae.2022.08.002

Rothermel, R. C. (1972). *A mathematical model for predicting fire spread in wildland fuels* (Vol. 115). Intermountain Forest & Range Experiment Station, Forest Service, US ….

Seidl, R., Schelhaas, M.-J., Rammer, W., & Verkerk, P. J. (2014). Increasing forest disturbances in Europe and their impact on carbon storage. *Nature Climate Change*, *4*(9), 806‑ 810. https://doi.org/10.1038/nclimate2318

Senande-Rivera, M., Insua-Costa, D., & Miguez-Macho, G. (2022). Spatial and temporal expansion of global wildland fire activity in response to climate change. *Nature Communications*, *13*(1), 1208.

Senf, C., Seidl, R., Knoke, T., & Jucker, T. (2024). *Taylor’s law predicts unprecedented pulses of forest disturbance under global change*.

Silva, J., Marques, J., Gonçalves, I., Brito, R., Teixeira, S., Teixeira, J., & Alvelos, F. (2022). A systematic review and bibliometric analysis of wildland fire behavior modeling. *Fluids*, *7*(12), 374.

Singh, R. D., Gumber, S., Joshi, H., & Singh, S. P. (2022). Allocation to tree bark in pine and oak species in fire affected mixed forests across the Northern Hemisphere. *Forest Ecology and Management*, *509*, 120081.

Spies, T. A., White, E., Ager, A., Kline, J. D., Bolte, J. P., Platt, E. K., Olsen, K. A., Pabst, R. J., Barros, A. M., Bailey, J. D., & others. (2017). Using an agent-based model to examine forest management outcomes in a fire-prone landscape in Oregon, USA. *Ecology and Society*, *22*(1).

Tapias, R., Climent, J., Pardos, J. A., & Gil, L. (2004). Life histories of Mediterranean pines. *Plant ecology*, *171*, 53‐ 68.

Taylor, S. W., Woolford, D. G., Dean, C. B., & Martell, D. L. (2013). Wildfire prediction to inform fire management : Statistical science challenges. *Journal of Environmental Sciences*.

Thonicke, K., Spessa, A., Prentice, I. C., Harrison, S. P., Dong, L., & Carmona-Moreno, C. (2010). The influence of vegetation, fire spread and fire behaviour on biomass burning and trace gas emissions : Results from a process-based model. *Biogeosciences*, *7*(6), 1991‐ 2011.

Varner, J. M., Shearman, T. M., Kane, J. M., Banwell, E. M., Jules, E. S., & Stambaugh, M. C. (2022). Understanding flammability and bark thickness in the genus Pinus using a phylogenetic approach. *Scientific Reports*, *12*(1), 7384.

Wagner, C. van. (1973). Height of crown scorch in forest fires. *Canadian journal of forest research*, *3*(3), 373‐ 378.

Zscheischler, J., Westra, S., Van Den Hurk, B. J., Seneviratne, S. I., Ward, P. J., Pitman, A., & Zhang, X. (2018). Future climate risk from compound events. *Nature Climate Change*, *8*(6), 469‐ 477.

# Supplementary Material

# Integrating Wildfire Risk into Robust Forest Management Optimization

Author Affiliation
Contact information
Clémence Labarre* (corresponding author)

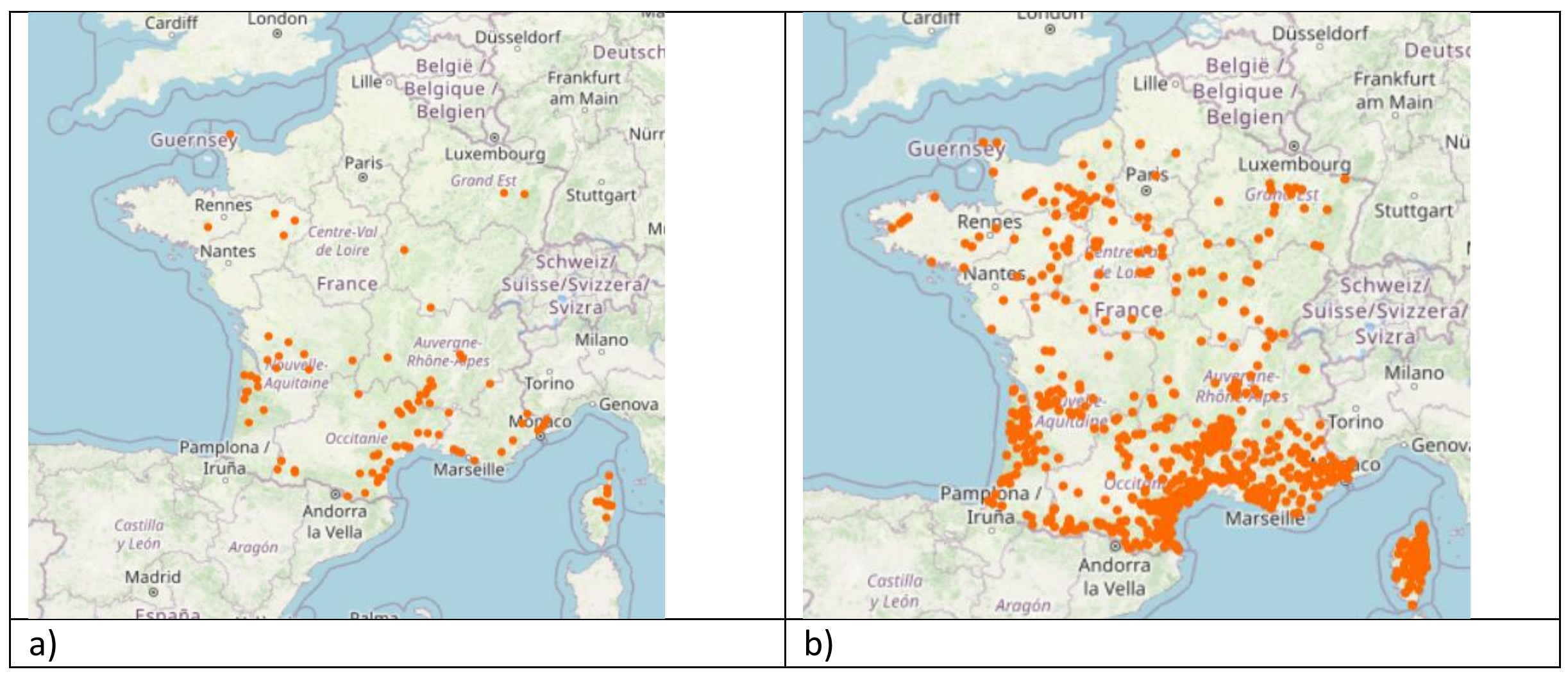

*Figure S1 Expansion of sever fire events between 2010 (A) and 2022 (B) Fire event burned an area greater than 20 ha (BDIFF database)*

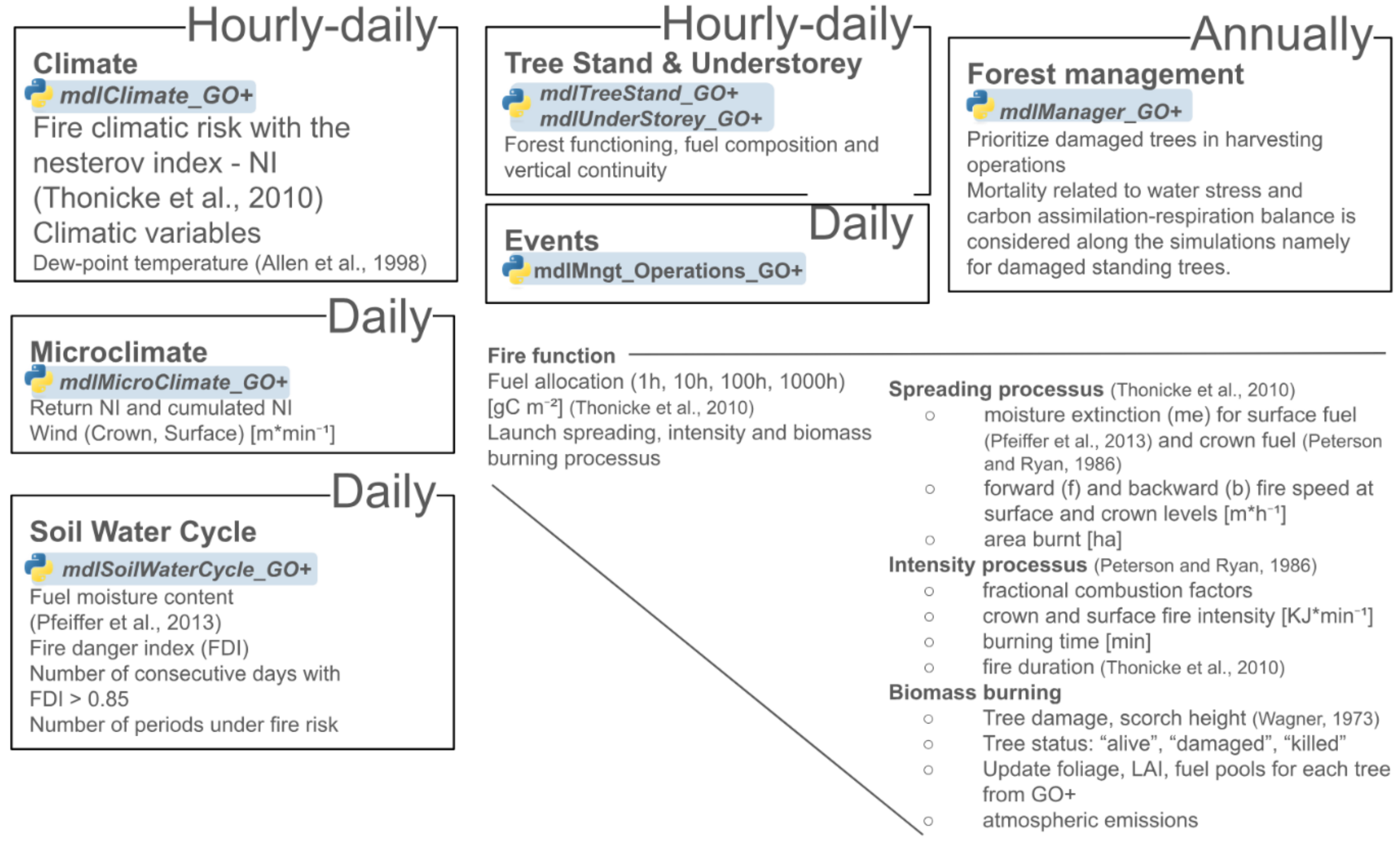

*Figure S2 Integration of fire related processes into the GO+ submodels. Equations describing the spread of fire in a forest stand have been integrated into the GO+ code (Peterson & Ryan, 1986; Pfeiffer, 2013; Rothermel, 1972; Wagner, 1973). They consider the forest ecosystem’s structure, biomass, ground fuel and water content. The weather-risk of fire is first evaluated at the hour and day scale, with climatic and water cycle parameters in the Climate and Soil water cycle modules (on the left). Then, the exposure of the forest and understorey is calculated in the Tree Stand and Understorey modules. The Operations module launch the fire event describing the spreading, intensity and burned biomass functions. The fire event triggers feedback effects on carbon cycle processes described in the Tree Stand module (Moreaux et al., 2020)}, including photosynthesis, respiration, litterfall, tree growth, soil carbon decomposition, canopy structure (e.g., height, crown base height, leaf area index), and biomass. Finally, the manager module evaluates tree damage and mortality before prioritizing salvaged trees in subsequent harvesting operations. This might influence in the Event module regeneration, thinnings, harvest, clearcutting, vegetation control, and mortality. References are given for the various added elements.*

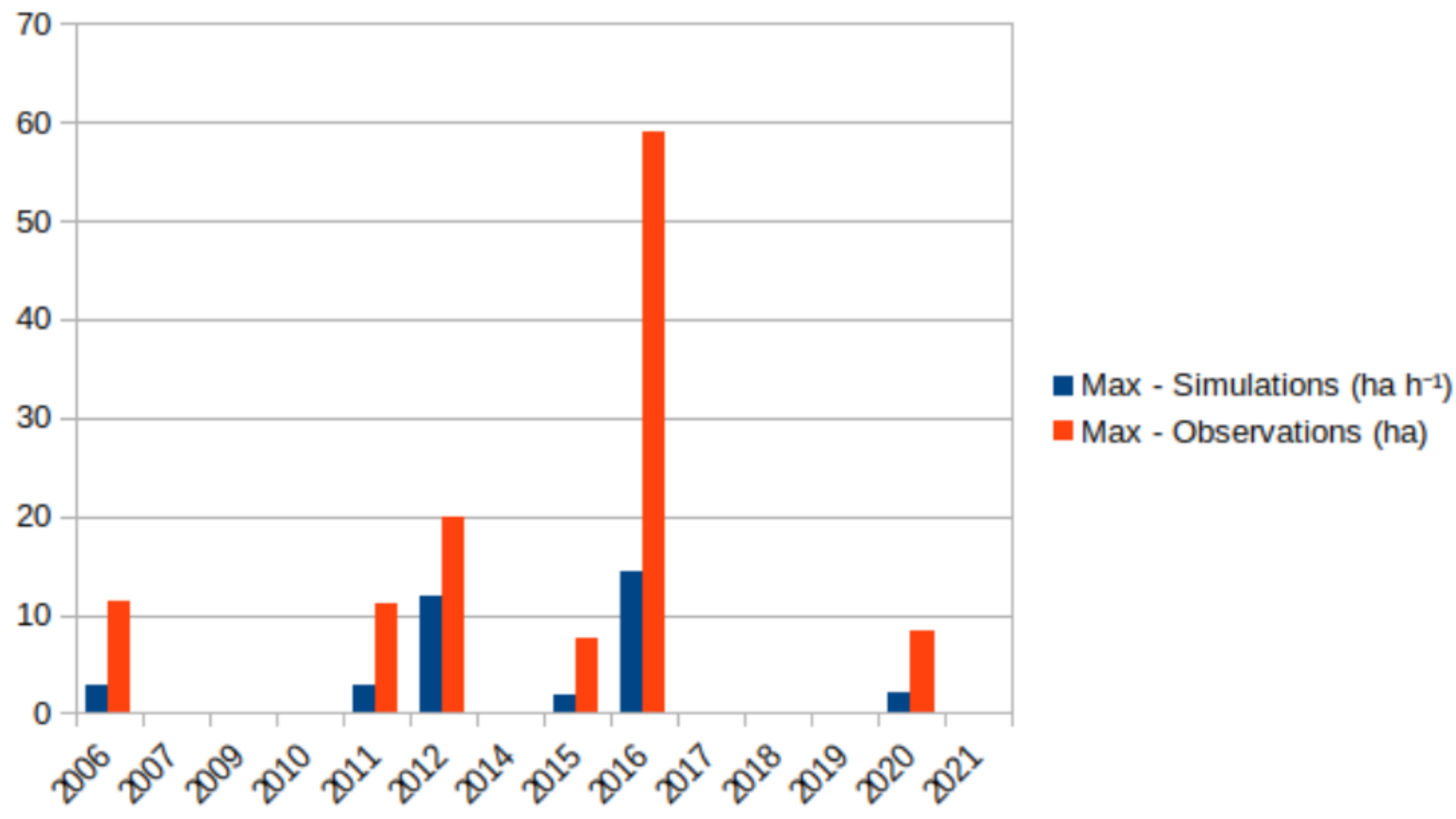

*Figure S3 Distribution of annual maximal historic and simulated fire events. (BDIFF database). The observed data, sourced from the BDIFF database, records fire events in the Landes de Gascogne region between 2004 and 2021, including fire dates, times, and extensions. Simulations were performed using the GO+Fire model for a 20-year-old maritime pine stand over a one-hour period, while observed data report total burned area. Given the dense firefighter network and rapid intervention, it is assumed that most observed fires lasted less than one hour, which justifies comparing observations and model predictions despite differing units (e.g.,ha vs ha $h^{-1}$).*

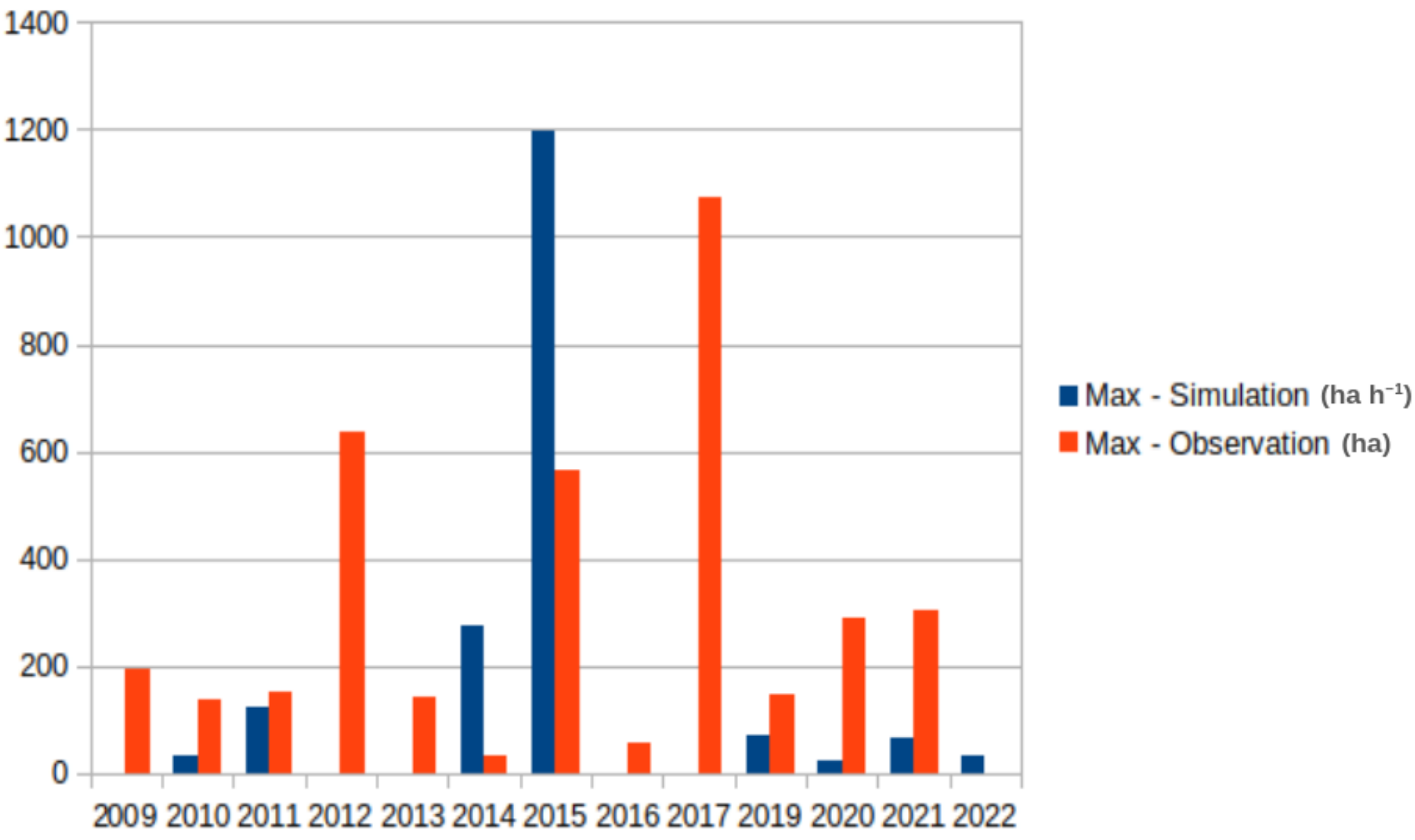

*Figure S4 Annual Maximum Area Burned (ha $h^{-1}$) of simulated and observed extreme events (>20.5 ha of burnt forest). The observed data, sourced from the BDIFF database, records fire events in the Landes de Gascogne region between 2004 and 2021, including fire dates, times, and extensions. The figures compared observed and simulated datasets focusing on fire events exceeding 20.5 ha. Simulations were performed using the GO+Fire model for a 20-year-old maritime pine stand over*

*a one-hour period, while observed data report total burned area. Given the dense firefighter network and rapid intervention, it is assumed that most observed fires lasted less than one hour, which justifies comparing observations and model predictions despite differing units (e.g.,ha vs ha $h^{-1}$).*

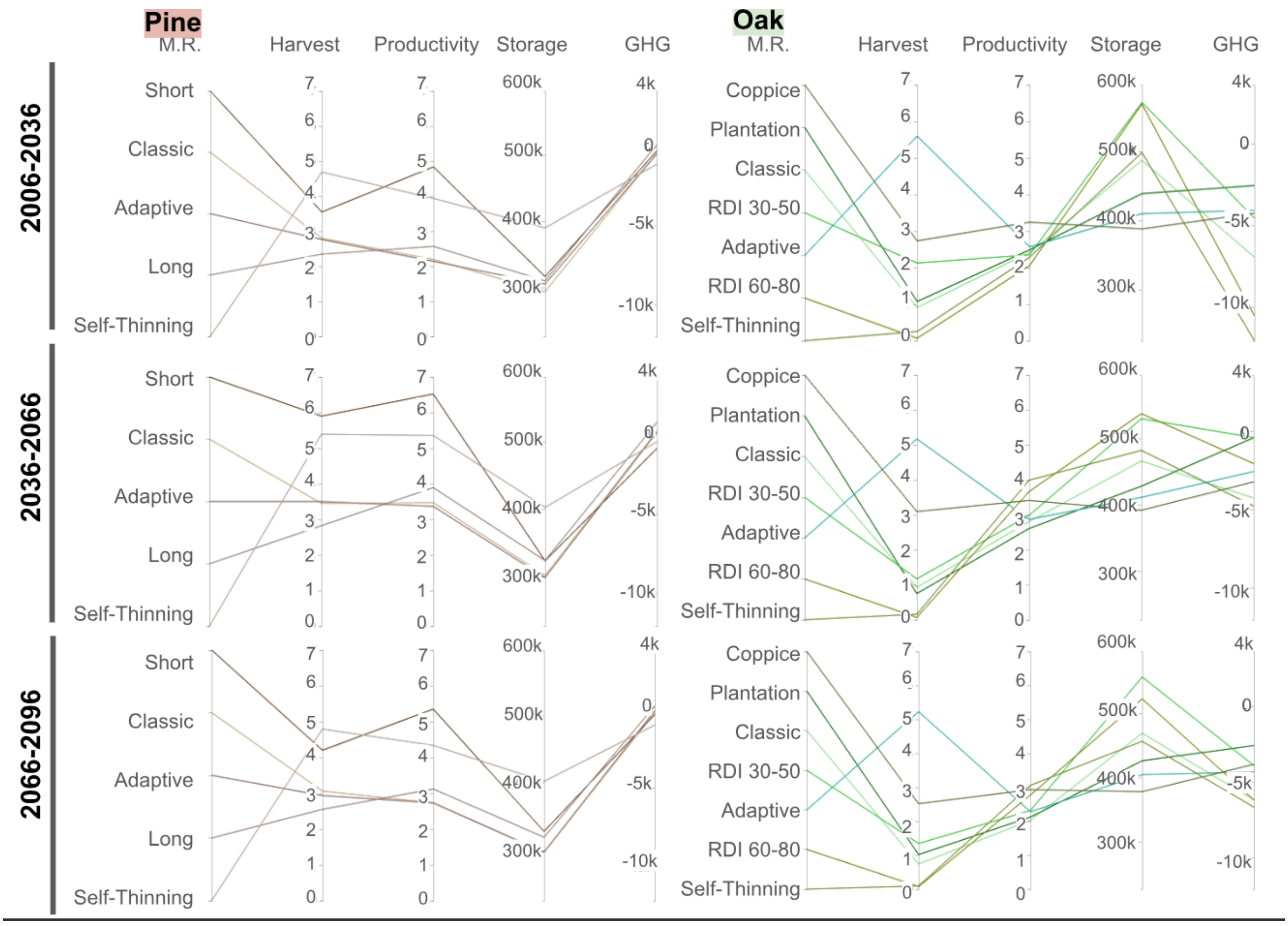

*Figure S5 Performance profile of forest management regimes under RCP 4.5.*
This parallel coordinate plot illustrates the performance of different management regimes across key indicators: timber production and forest productivity ($m^3$ $ha^{-1}$ $yr^{-1}$), as well as carbon storage and ecosystem greenhouse gas budget (gC-$CO_{2e}$ $m^{-2}$ over 30 years). Each management regime was simulated over two sites with contrasting rainfall regime, and their climate responses were averaged. Results are presented for three climatic periods: 2006–2036 (A), 2036–2066 (B), and 2066–2096 (C); and two species: *Pinus pinaster* and *Quercus robur.*

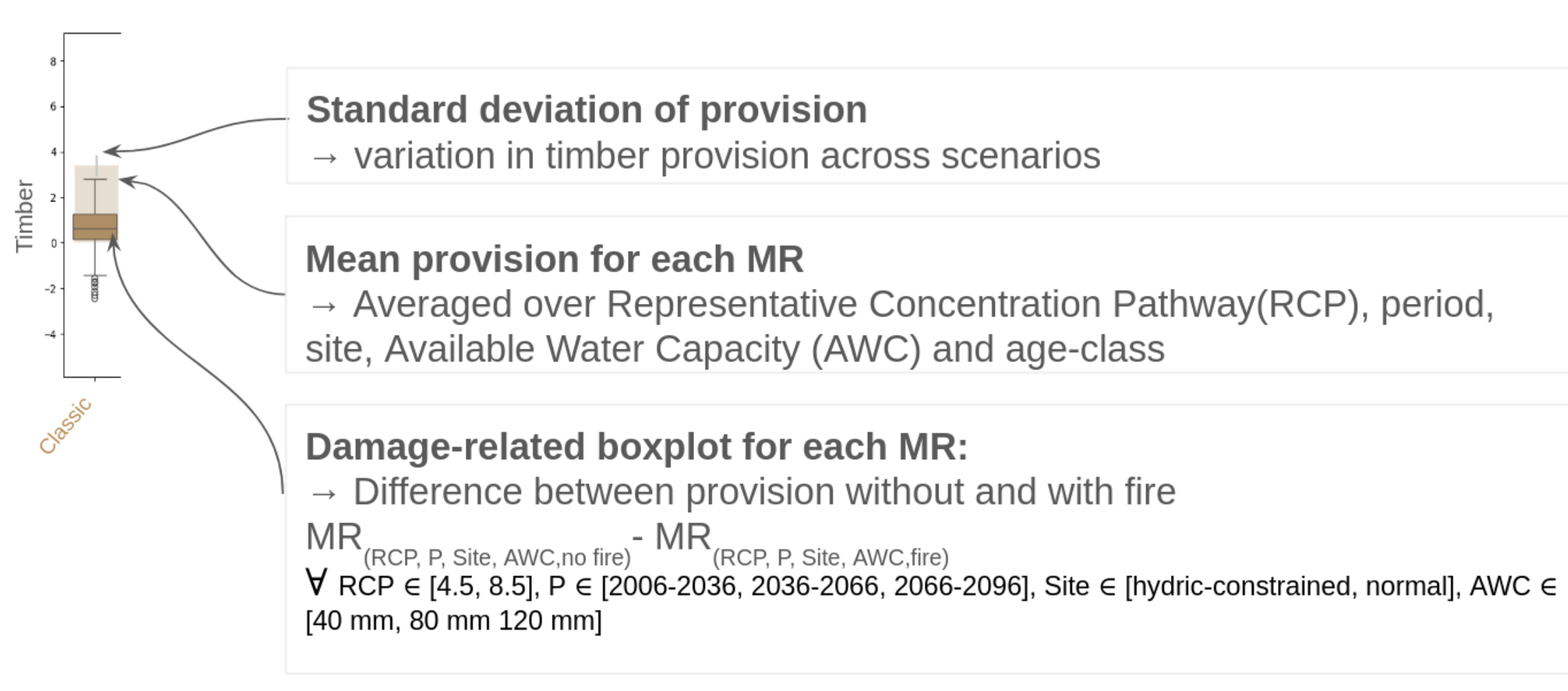

*Figure S6 Explanatory Key for Figure 4*

*Table S1: Characteristics of forest management regimes*

The forest management regimes are described as itineraries starting from regeneration and running until the next final clearcut thus covering the entire life cycle of the tree stand. The tree stand are monospecific and even aged. Throughout the life of a stand, tree density is thus controlled by regeneration, mortality, thinnings and final cutting. The latter are triggered according to rules summarized in the table below.

A - Maritime pine – *Pinus pinaster*

| Regime | Soil preparation (fraction affected , soil depth (m)) (1) | Planted density ($ha^{-1}$) | Spacing (dominant height, final density) (year) | Fraction removed and circumference threshold at thinnings (%) | Circumference threshold at final harvest (cm) | Parts of trees harves |
|---|---|---|---|---|---|---|
| Classic | 1, 0.2 | 2500 | 5, 10 | 50% at 7-year, then 30% at 50, 70, 85, 100 cm | 120 | Trunk |
| Short | 1, 0.2 | 2500 | 5, 10 | 30% at 50, 70 cm | 85 | Trunk, branch stump |
| Long | None | 1250 | 5, 10 | 30% at 55, 75, 95, 115 cm | 150 | Trunk |
| Adaptive | 1, 0.2 | 1250 | 5, 10 | 30% when cumulated stress index (2) ≥ 0.5 or circumference at 50, 70, 85, 100 cm | 100 | Trunk |
| Self-thinning | None | 3500 | Self-thinning | Self-thinning (3) | – | Trunk |

(1) The soil is prepared by mechanical disking.
(2) The stress index used for the Adaptive regimes is given by eq. 29 in Moreaux et al. (2020) as the complement to unity of the ratio of actual to potential annual tree stand transpiration.
(3) The applied natural mortality rate is based on forest inventories. For all management regimes, this mortality rate is 0.5% yr-1 if the mean tree diameter is less than 0.55 m, and 0.02% yr-1 otherwise.

B – Sessile and pedunculate oaks – *Quercus petraea and robur*

Here, depending on the specific regime, the process of thinning may be initiated according to age, diameter, height, competition (RDI) or stress index thresholds. RDI, the ratio of current density to the maximal sustainable density at which mortality occurs, is a metric employed in the Adaptive and RDI regimes.

| Regime | Soil preparation (fraction affected, soil depth (m)) (1) | Regeneration type (trees per ha) | Spacing (dominant height, final density) | Thinning rule | Fraction removed by thinnings | Rotation period | Parts of trees harvested |
|---|---|---|---|---|---|---|---|
| Classic | 1, 0.1 m | Natural (5000) | 6 m (3000) | Height = 14 m, then every 10 to 16 years | 66% then 30–40% | 170–180 years | Trunk |
| Adaptive | 1, 0.1 m | Natural (5000) | 6 m (3000) | RDI < 0.65 or Stress index cumulated over 7 years > 3.5 | to RDI = 0.40, | 170–180 years | Trunk |
| Plantation | 1, 0.1 m | Plantation (1600) | 12 m (500) | every 8 to 10 years | 20 to 40% | 90 years | Trunk |
| Coppice | 1, 0.1 m | Re-sprouting (3500) | – | age = 10-year | 64% | 45 years | Trunk |
| RDI(2) 30–50 | 1, 0.1 m | Natural (3000) | 6 m (1000) | RDI > 0.5 | to RDI = 0.3 | 150–180 years | Trunk |
| RDI(2) 60–80 | 1, 0.1 m | Natural (3000) | 6 m (900) | RDI > 0.8 | to RDI = 0.6 | 150–180 years | Trunk |
| Self-thinning | 1, 0.1 m | Natural (5000) | – | Self-thinning(3) | variable | 200 years | Trunk |

(1) The “soil preparation” corresponds to the disturbance made by logging machinery.
(2) Relative Density Index
(3) In addition to the regular thinnings, self mortality is implemented in all regimes according to Ningre et al. (Ann For Sci 2019Ningre, F., Ottorini, J. M., & Le Goff, N. (2019). Size-density trajectories for even-aged sessile oak (Quercus petraea (Matt.) Liebl.) and common beech (Fagus sylvatica L.) stands revealing similarities and differences in the mortality process. Annals of Forest Science, 76, 1-20.)

## Appendix 1. Document on the Integration of Fire Mechanisms into GO+, a Process-Based Model of Forest Growth and Management

Depending on the thinning regime, vegetation control, soil preparation and harvest fraction, management alternatives could result in different amounts of dead and live fuels and their water content, which are the main determinants of stand fire susceptibility. Here, our goal was to investigate how fire damages to ecosystem services could be affected by forest management alternatives across a 1 million hectares forest showing substantial variations in soil available water content (AWC) and climate. To this end, we derived a fire damage submodel from simple, widely used fire models (Pfeiffer, 2013; Rothermel, 1972; Thonicke et al., 2010). More sophisticated forest fire models based on the 3D canopy structure would have required a realistic 3D description of the canopy, which is not available from the GO+ homogeneous layer approach, and would have raised intractable problems for implementing multidecadal simulations of different management alternatives under climate scenarios.

The SPITFIRE model (Thonicke et al., 2010) simulates forest fire intensity, duration and expansion on a daily basis and was devoted to the coupling with DGVM global models such as LPJ, JULES or ORCHIDEE models. The main assumptions of SPITFIRE concern the fire duration that is restricted to a maximum of 6 hours and its limitation to surface fire. They were partially relaxed in the LPJ-LMfire version (Pfeiffer, 2013). Both models describe the state of the fuels in a given PFT, the ignition probability and estimates the fire danger based on the equations and model proposed by Rothermel (Rothermel, 1972). The Rothermel model simulates the intensity and duration of fires as well as their impact on the fuel beds and tree mortality and damage. The fuel is partitioned among the ground floor and the trees and classified in size classes according to the time needed for reaching its equilibrium moisture, 1, 10, 100 or 1000 hours.  Here, we coupled a simple fire model derived for Rothermel (Rothermel, 1972) with the process-based model of forest growth GO+.

The GO+ model described the forest as a two-layer vegetation canopy, i.e. the tree layer and ground vegetation. The ROTH-C and ECOSS models describe the carbon and nitrogen cycles in soils and are integrated in GO+ v3.2. They simulate the dynamics of the decomposable plant material and resistant plant material deposited into the ground y vegetation layers.  Go+ solved the processes of radiation and heat exchange, evapotranspiration, rainfall interception, photosynthesis and respiration for each canopy layer at an hourly time step. It simulates the phenology and growth of trees and vegetation based on the carbon, water and nitrogen availability. It has been calibrated for a number of tree species (*Pinus, Pseudotsuga, Quercus, Fagus, Eucalyptus*) and ground vegetation types (grasses, woody plants, legumes). The French Forêts-21 project (https://forets21.inra.fr/pelican3.1/modele.html) has made it possible to include the main management operations (plantation, regeneration, soil preparation, spacing, vegetation control, thinning and clearcutting) and to build up from them 30 management alternatives adapted to French forests. We included the basic equations of the Rothermel model in the GO+ model called “Fire” and adapted them using forest specific equations mainly taken from Peterson and Ryan (Peterson & Ryan, 1986).

- The stocks of dead fuel are calculated at a daily resolution by the ROTH-C model included in GO+. The 1h fuel component corresponds to the aboveground Decomposable Plant Material and the 10h 100h and 1000h fuels components summed to the aboveground Resistant Plant

Material component. The moisture of each fuel fraction is calculated hourly based upon the water balance model proposed by Pfeiffer et al. (Pfeiffer, 2013).

- The stock of the live fuels is also provided by GO+ from the trees and understorey sub-models. The live fuels stock is updated at a daily resolution for the tree foliage and understorey or annually for the tree stems and branches respectively. The moisture content of both understorey and trees parts are estimated from their leaf water potential calculated hourly by GO+. The basic density of live fuels is calculated using the size and biomass amount accumulated in each canopy layer.
- Damage is calculated for each individual tree for each fire. Damages inflicted on an individual tree are two-fold: the crown can be partially or totally burned and the cambial tissues of the trunk are damaged by excessive heating and burning. We used the Peterson and Ryan equation (8) that combines both effects into a single equation. Accordingly, the probability of killing a tree depends on the relative scorch height and the fire duration relative to the critical time needed for destroying the whole bark and cambium. When not killed, a burned tree may survive if the majority of its crown is left intact and the fire duration is less than the critical time required to burn the bark and cambium. Damage of surviving trees is reducing their leaf area and degrades the quality of their stemwood.
- Primary production, canopy respiration, nitrogen uptake, and , in turn, the growth of surviving trees are therefore impacted by fire. Dead, unburned biomass fractions from surviving or killed trees are deposited in the soil compartment as aboveground litter or may be harvested. In the latter case, the harvested material is assumed to be degraded and therefore can only be subsequently used for pulp or energy. The trunks of damaged trees are preferably felled during the next thinning, starting with most damaged trunks. If the fire burns all the trees in the stand, a new stand is regenerated according to the specifications of the management scenarios and after one year of logging and soil preparation.

Though GO+ runs with hourly resolution, the fire is simulated as a single, instantaneous, event. In other words, we did not attempt at simulating the fire dynamics hour after hour. The effect of fire is calculated as the difference between unburnt and burnt forest stand times the total area affected by fires. The fire effect is thus quantified on a number of ecosystem services and production through the corresponding variables:

- Wood production and its distribution among product categories
- Carbon stocks in soil, biomass, and harvested products
- Net balance of grenhouse gases exchange
- Energy fluxes
- Water balance components: evapotranspiration and deep runoff

Here, only the production, carbon stocks and GHG fluxes were analysed.

## Evaluation.

We used the recording of the SDIS Gironde for evaluation how well the model reproduces the forest fire observed in the Gironde department from 2004 to 2021. The variables available from SDIS are the date and time of fires and their extension. We simulated 40 fires recorded from 2005 to 2020 and their impact on a typical 20-year old maritime pine stand. We compared the surface burnt during one hour between the values simulated by GO+ and the observed values. The total area burnt was not used in this comparison since GO+F does not include the extinction of fire by firemen which is the case for a large majority of fires in SouthWestern France. GO+ limits the duration of fire to 4 hours. Instead, the area burnt per hour were calculated and compare with observations.

### Experiment.

In order to quantify the damages caused by forest fires to forest ecosystem services, we have first estimated the damages at the plot scale and then upscaled the damaged for the entire forest area.

At plot scale, the damages are quantified as the difference between a "control" and a "fire" runs of GO+. They were implemented for 3 periods extending respectively from 2006 to 2036, 2037-2066 and 2067-2096. For each period the simulations were initialized for 6 to 10 age classes depending on the management considered and subsequently averaged over age classes. Two grid points, corresponding to the driest and the wettest areas of the forest, were implemented using the medium class of soil WHC, that is 80 mm. Two RCP scenarios were used, RCP45 and RCP85.

One fire was simulated for each 30-year period. For establishing the fire date and time in each period, a first run of a standard, mature stand of trees was performed to calculate the time course of the fire danger index from 2006 to 2100. ThenAfter a fire, the whole tree stand is logged when more than 40% of the stems are killed. The trunks are harvested and a new stand is regenerated after a delay of one year. Otherwise, the tree stand is maintained with a reduced leaf area depending on the number of surviving trees and the proportion of unburned foliage calculated for each surviving stem. In both cases, the logs harvested from burned trees are considered as degraded and used only for pulp or energy.

**Table S2: Ignition dates**

| Period | Starting in | Grid point | RCP | Year of fire | Day of fire |
|---|---|---|---|---|---|
| period 1 | 2006 | 7664 | 4.5 | 2021 | 228 |
| period 2 | 2036 | 7664 | 4.5 | 2051 | 219 |
| period 3 | 2066 | 7664 | 4.5 | 2081 | 215 |
| period 1 | 2006 | 7664 | 8.5 | 2018 | 246 |
| period 2 | 2036 | 7664 | 8.5 | 2054 | 203 |
| period 3 | 2066 | 7664 | 8.5 | 2078 | 231 |
| period 1 | 2006 | 8357 | 4.5 | 2018 | 250 |
| period 2 | 2036 | 8357 | 4.5 | 2050 | 226 |
| period 3 | 2066 | 8357 | 4.5 | 2081 | 227 |

### Damage index.

The damage index, $I_p$, is calculated as;

$l_p = D \times A_b \times f_{fire}$ Where $D$ is the penalty due to fire for a given variable and is expressed per unit area, $A_b$ the area burned and $f_{fire}$ the frequency of fire danger calculated as the number of days with FDI>0.85 over the total of days of the period (10957 days);

the damage can be either the loss of timber production (10^3m3), biomass production (ton C ), water deep runoff (m3 H2O), carbon storage in biomass and soil (ton C ).

**Table S3: Fire event characteristics from 2006 to 2100 simulation**

| Year | Species | Fire # | Duration | DOY | Fuel1h (gC $m^{-2}$) | |
|---|---|---|---|---|---|---|
| 2006 | Pinaster | 1 | 5 | 232 | 300.16198 | |
| 2006 | Pinaster | 2 | 15 | 240 | 337.39485 / 7 | |
| 2007 | Pinaster | 3 | 3 | 239 | 221.36675 | |
| 2008 | Pinaster | 4 | 7 | 238 | 249.82950 / 3 | |
| 2008 | Pinaster | 5 | 20 | 246 | 278.64659 / 293.50803 | |
| 2008 | Pinaster | 6 | 4 | 275 | 249.60728 / 4 | |
| 2009 | Pinaster | 7 | 3 | 273 | 313.58287 / 9 | |
| 2009 | Pinaster | 8 | 16 | 280 | 208.78797 / 2 | |
| 2010 | Pinaster | 9 | 1 | 220 | 3 | |
| … | … | … | … | … | … | |
| 2099 | Pinaster | 379 | 3 | 264 | 111.84618 / 112.52605 | |
| 2099 | Pinaster | 380 | 1 | 271 | 4 | |
| 2100 | Pinaster | 381 | 20 | 219 | 136.94185 / 9 | |
| 2100 | Pinaster | 382 | 2 | 244 | 138.17934 / 8 | |
| 2100 | Pinaster | 383 | 9 | 247 | 147.33592 | |
| 2100 | Pinaster | 384 | 8 | 257 | 154.86807 | |
| 2100 | Pinaster | 385 | 3 | 280 | 170.57938 / 7 | |
| 2100 | Pinaster | 386 | 9 | 284 | 183.91147 / 5 | |

*Table S4: Multivariate Random Forest Model Performance and Variable Importance for Fire Response Metrics. We used a Multi-output Random Forest model (via MultiOutputRegressor with RandomForestRegressor as the base estimator) to simultaneously predict changes in targets including damage on harvest, GHG budget, carbon storage, productivity, as well as fire probability from a shared set of explanatory variables: climate change scenario, 30-year climatic period, management regime, species, available water content, site, year.*

A. Model performance Metrics

| Target Variable | $R^2$ | MAE | RMSE |
|---|---|---|---|
| $\Delta_{(}$Harvest$_{)}$ | 0.817 | 0.315 | 0.520 |
| $\Delta_{(}$GHG Budget$_{)}$ | 0.527 | 103.495 | 194.173 |
| $\Delta_{(}$Carbon Storage$_{)}$ | 0.759 | 624.035 | 926.126 |
| $\Delta_{(}$Productivity$_{)}$ | 0.587 | 0.485 | 0.713 |
| Fire Probability ($P_{(}$fire$_{)}$) | 0.999 | 2.377 | 10.800 |

B. Average Feature Importance Across All Targets

| Variable | Importance |
|---|---|
| Management Regime | 0.361749 |
| Year | 0.315781 |
| Site | 0.136461 |
| 30-year climatic period | 0.059648 |
| Species | 0.049006 |
| Available Water Capacity | 0.044629 |
| RCP | 0.032276 |

# References

Senf, C., Seidl, R., Knoke, T., Jucker, T., 2025. Taylor's law predicts unprecedented pulses of forest disturbance under global change. Nat Commun 16, 6133. https://doi.org/10.1038/s41467-025-61585-5.